\documentclass[12pt]{article}
\usepackage{amsfonts, latexsym,amssymb,amsmath,amstext,epsfig}

\title{Uniform convergence of hypergeometric series}
\author{Raimundas Vid\=unas\footnote{Primarily supported by NWO,
project number 613-06-565. Also supported by the ESF NOG project, and the 21
Century COE Programme "Development of Dynamic Mathematics with High
Functionality" of the Ministry of Education, Culture, Sports, Science and
Technology of Japan. }}

\catcode`\@=11 \@addtoreset{equation}{section} \catcode`\@=12

\newtheorem{theorem}{Theorem}[section]
\newtheorem{lemma}[theorem]{Lemma}

\newtheorem{conclude}[theorem]{Corollary}

\newtheorem{example}[theorem]{Example}
\newcommand{\proof}{{\bf Proof. }}
\newcommand{\QED}{\hfill \mbox{QED.}\\}
\renewcommand{\labelenumi}{(\roman{enumi})}
\newcommand{\refi}[1]{\mbox{\em (#1)}}

\newcommand{\hpg}[5]{{}_{#1}\mbox{\rm F}_{#2}\!
  \left(\left.{#3 \atop #4}\,\right|#5 \right) }
\newcommand{\hpgo}[2]{{}_{#1}\mbox{\rm F}_{\!#2}}

\newcommand{\equal}{\!\!=\!\!}
\newcommand{\wtilde}{\widetilde}
\newcommand{\wideha}{}

\newcommand{\re}[1]{\mbox{\rm Re}\,#1}
\newcommand{\ree}[1]{\mbox{\rm Re}\left(#1\right)}
\newcommand{\scre}[1]{\mbox{\scriptsize\rm Re}\left(#1\right)}
\newcommand{\CC}{\mathbb C}
\newcommand{\ZZ}{\mathbb Z}
\newcommand{\RR}{\mathbb R}
\newcommand{\Zplus}{\mbox{\bf Z}_+}
\newcommand{\bS}{{\cal U}}

\date{}

\begin{document}

\maketitle

\begin{abstract}
The considered problem is uniform convergence of sequences of hypergeometric
series. We give necessary and sufficient conditions for uniformly dominated
convergence of infinite 
sums of proper bivariate hypergeometric terms.  
These conditions can be checked algorithmically. Hence the results
can be applied in Zeilberger type algorithms for
nonterminating hypergeometric series.
\end{abstract}

\section{Introduction}

In this paper we study uniform convergence of sequences of hypergeometric
series. We consider the sequences $\bS(n)=\sum_{k=0}^{\infty} u(n,k)$ of
hypergeometric series such that $u(n,k)$ is a proper hypergeometric term in
$n,k$. We assume that the individual series $\bS(n)$ are nonterminating for
large enough $n$. The underlying field is the complex numbers.

Recall \cite{abrpet02} that a bivariate sequence $u(n,k)$ is a {\em
hypergeometric term} if both quotients $u(n+1,k)/u(n,k)$ and
$u(n,k+1)/u(n,k)$ can be realized as rational functions of $n,k$.

A bivariate sequence $u(n,k)$ is a {\em proper term} if there exist: a
non-negative integer $p$; complex constants $\xi,\theta$; $b_1,\ldots,b_p$;
integers $\alpha_1,\ldots,\alpha_p$; $\beta_1,\ldots,\beta_p$; and a
polynomial $P(n,k)$ such that
\begin{equation} \label{eq:properterm}
u(n,k)=P(n,k)\,\frac{\xi^n\theta^k}{k!}\,(b_1)_{\alpha_1n+\beta_1k}\,\cdots\,(b_p)_{\alpha_p\,n+\beta_p\,k},
\end{equation}
where $(a)_m$ is the Pochhammer symbol:
\begin{equation}
(a)_m=\left\{ \begin{array}{rl} a\,(a+1)\cdots(a+m-1), & \mbox{if } m>0,\\
1, & \mbox{if } m=0,\\ 1/(a-1)\,\cdots\,(a-|m|), & \mbox{if } m<0.
\end{array}\right.
\end{equation}
In general, $(a)_m=\Gamma(a+m)/\Gamma(a)$. To avoid consideration of undefined
or terminating series, we assume that:
\begin{itemize}
\item[\it(I)] For those $j$ with $\alpha_j>0$ or $\beta_j>0$ we have
$b_j\not\in\{0,-1,-2,\ldots\}$.
 \item[\it(II)] For those $j$ with $\alpha_j<0$ or $\beta_j<0$ we have
$b_j\not\in\{0,1,2,\ldots\}$.
\end{itemize}
In particular, for those $j$ with $\alpha_j\beta_j<0$ we must have $b_j\not
\in\ZZ$.
More generally, one may consider proper hypergeometric terms with the
Pochhammer symbols in (\ref{eq:properterm}) replaced by other factorial-type
functions, such as the gamma function of linear arguments in $n,k$ with
integer coefficients to $n,k$, or the {\em nonvanishing rising factorial}
\cite{abrpet02}.

A bivariate sequence $u(n,k)$ is {\em holonomic} if the generating function
$\sum_{n,k\ge 0} u(n,k)\,x^ny^k$ and all its partial derivatives generate a
finite-dimensional vector space over the field of rational functions in
$x,y$.

Proper terms are hypergeometric and holonomic \cite[Theorem 3]{abrpet02}.
Since our coefficient field is algebraically closed, 
any holonomic hypergeometric term is {\em conjugate} to a proper term
\cite[Theorem 14]{abrpet02}, \cite{hou}. This means that there are polynomials
$f_1,f_2,g_1,g_2$ in $n,k$, not identically zero, and a proper term
$\widetilde{u}(n,k)$, such that
\begin{equation}
f_1(n,k)\,u(n+1,k)=f_2(n,k)\,u(n,k), \qquad
g_1(n,k)\,u(n,k+1)=g_2(n,k)\,u(n,k),
\end{equation}
and these identities hold with $u$ replaced by $\widetilde{u}$ as well.


If a term $u(n,k)$ is holonomic, then it satisfies difference equations (in
one or both variables) whose coefficients are dependent only on $n$
\cite[Chapter 4]{aeqb}. If the term $u(n,k)$ is proper, and for any $n$ the
sum $\bS(n)=\sum_{k=0}^{\infty} u(n,k)$ is terminating, Zeilberger's
algorithm gives a recurrence relation with respect to $n$ for $\bS(n)$. The
crucial step in Zeilberger's algorithm is to derive a recurrence relation
\begin{equation} \label{eq:zeilb}
L(n)\,u(n,k)=R(n,k+1)-R(n,k),
\end{equation}
where $L(n)$ is a linear difference operator with coefficients in $n$ only,
and $R(n,k)$ is a hypergeometric term. The linear recurrence is derived by
summing (\ref{eq:zeilb}) over all $k$; the right hand-side simplifies due to
telescoping summation.

When generalizing Zeilberger's algorithm to nonterminating hypergeometric
series, one needs to make sure that the series $\bS(n)=\sum_{k=0}^{\infty}
u(n,k)$ converges uniformly, so to justify manipulation of 
(\ref{eq:zeilb}). This paper gives criteria to decide uniform convergence of
$\bS(n)$. In \cite{hpgzeilberg} these criteria are used for the Zeilberger
type algorithms for nonterminating hypergeometric series. Generalization of
$q$-Zeilberger algorithm to nonterminating basic hypergeometric series is
considered in  \cite{chenhoumu}.

The main result of this paper is the sufficient and necessary conditions for
uniformly dominated convergence of sequences $\bS(n)=\sum_{k=0}^{\infty}
u(n,k)$ of infinite sums 
of proper hypergeometric terms $u(n,k)$. Uniformly dominated convergence
is defined by the Weierstrass M-test; see Lemma \ref{weierstrassm} below.

The main result is presented in Section \ref{sec:theresult},
and proved in Section \ref{mainproof}. In Section \ref{examplez}, examples of
application of the main result are given.  In Section \ref{sec:basic},
the crucial technical Lemma \ref{hpguc} and a few asymptotic expressions
for the gamma function are presented. In Section \ref{sec:otherprel} we provide several
other intermediate results. In Sections \ref{sec:notation} and
\ref{sec:enotation} we specify the form of hypergeometric series under
consideration, and define the notation we use. In Section \ref{gtsection} we
consider the least straightforward (from computational point of view) part
of the main theorem more closely. 


\section{Basic preliminary results}
\label{sec:basic}

Throughout the paper, let $\Zplus$ denote the set of non-negative integers.
We make the convention that $0^0=1$, which is the proper continuous limit of
the function $|x|^x$.

As the criterium for uniform convergence of function series,
we use the Weierstrass M-test formulated here below. (We apply it with
$E=\Zplus$.)
\begin{lemma} \label{weierstrassm}
Let $f_0(x),f_1(x),f_2(x),\ldots$ be a sequence of
complex-valued functions on a set $E$. 
If there exists a sequence $M_0,M_1,M_2,\ldots$ of real constants such that
$|f_j(x)|\le M_j$ for any $x\in E$ and all $j\in\Zplus$, and the series
$\sum_{j=0}^{\infty}M_j$ converges, then the function series
$\sum_{j=0}^{\infty}f_j(x)$ converges uniformly on $E$.
\end{lemma}
We refer to a function series that satisfies the sufficient condition of
this criterium as a {\em uniformly dominated convergent} series. In plain
terms, the condition is that the series is uniformly bounded (or majorized)
by an absolutely convergent series.

The following lemma gives us a strategy to determine uniformly dominated
convergence of sequences of nonterminating hypergeometric series.
\begin{lemma} \label{hpguc}
Let $u(n,k)$ denote a hypergeometric term in $n,k$. We assume that the
hypergeometric series $\bS(n)=\sum_{k=0}^{\infty} u(n,k)$ is nonterminating
for large enough $n$.
The series sequence $\bS(n)$ is uniformly dominated convergent
if and only if the following conditions hold:
\begin{itemize}
\item[(a)] For any $n\ge 0$, the series $\bS(n)$ converges absolutely.
 \item[(b)] The termwise limit $\sum_{k=0}^{\infty} \lim_{n\to\infty} u(n,k)$
exists and converges absolutely.
 \item[(c)] For any function $N:\Zplus\to\Zplus$ such that $N(k)\sim C_0\,k^p$
for some real $p>1$ and $C_0\neq 0$, the series $\sum_{k=0}^{\infty}
u(N(k),k)$ converges absolutely.
 \item[(d)] For any function $N:\Zplus\to\Zplus$ such
that $N(k)\sim C_0\,k^p$ for some real $p\in\!(0,1)$ and $C_0\neq 0$, the
series $\sum_{k=0}^{\infty} u(N(k),k)$ converges absolutely.
 \item[(e)] For any function $N:\Zplus\to\Zplus$ such that
$N(k)=\lambda\,k+\omega(k)$ for some real $\lambda>0$, with either
$\omega(k)=O(1)$ or $\omega(k)\sim C_0\,k^p$ for some real $p\in(0,1)$,
$C_0\neq 0$, 
the series $\sum_{k=0}^{\infty} u(N(k),k)$ converges absolutely.
\end{itemize}
\end{lemma}
\proof The conditions are necessary because a uniformly bounding series
would be a majorant for the indicated series as well. The limit in {\em (b)}
exists because, for each $k$, $|u(n,k)|$ is a monotonic function of $n$ for
$n$ large enough.

To prove the sufficiency, we may assume that for $k$ large enough $u(n,k)$
is not a constant in $n$. 
Let $z(k)=\sup_{n\ge 0} |u(n,k)|$. Then the series $\sum_{k=0}^{\infty}
z(k)$ is a precise uniform majorant for $\bS(n)$. The series sequence
$\bS(n)$ is uniformly dominated convergent if and only if the series
$\sum_{k=0}^{\infty} z(k)$ converges.

Let $h(\nu,\kappa)$ be the rational function of two complex variables 
equal to \mbox{$u(\nu+1,\kappa)/$}{$u(\nu,\kappa)$} for positive integer
values of $\nu$ and $\kappa$. Note that $h(\nu,\kappa)$ is nonzero and
well defined for positive integers $\kappa$ and large enough integers $\nu$,
because $u(n,k)=0$ would imply that the hypergeometric series are
terminating or undefined for large enough $n$. The function $h(\nu,\kappa)$
may be complex-valued, but the variables $\nu,\kappa$ are assumed to be
real.

For each non-negative integer $k$, we have that either
$z(k)=\lim_{n\to\infty} |u(n,k)|$, or
\[
\mbox{$z(k)=|u(n_0,k)|$ and $|h(n_0,k)|\le1$, $|h(n_0-1,k)|\ge1$ for some
integer $n_0$}.
\]
In the latter case, the rational function $|h(\nu,k)|^2$ of $\nu$ (with $k$
fixed) has either a pole on the interval $\nu\in[n_0\!-\!1,n_0]$, or it is
continuous and therefore achieves the value 1 on the same interval.

Let $\widehat{h}(\nu,\kappa)$ be the denominator of $|h(\nu,\kappa)|^2$.
Since the series $\bS(n)$ are not constant when $n$ large enough, we have
that $|h(\nu,k)|$ is not a constant function of $\nu$ for all large enough
$k$. Let $\nu_1(\kappa),\ldots,\nu_m(\kappa)$ be the positive real algebraic
functions, which are solutions of the algebraic equations
\mbox{$|h(\nu,\kappa)|=1$} or $\widehat{h}(\nu,\kappa)=0$, and are defined
for large enough $\kappa$. For $j=1,\ldots,m$, let $N_j(k)$ be the
integer-valued function
\[
N_j(k)=\left\{\begin{array}{rl} \lfloor \nu_j(k) \rfloor, &
\mbox{if } |h(\lfloor\nu_j(k)\rfloor,k)|\le 1, \\
\lceil \nu_j(k) \rceil, & \mbox{if } |h(\lfloor\nu_j(k)\rfloor,k)|>1.
\end{array}\right.
\]
All these functions satisfy the assumption of one of the last three
conditions. 
The functions $N_j(k)$ give candidates for ``local" maximums of $|u(n,k)|$,
as $n$ varies over the discrete set of positive integers and $k$ is fixed.

For large enough $k$, the candidates for $z(k)$ are $|u(0,k)|$,
$\lim_{n\to\infty} |u(n,k)|$, and $|u(N_j(k),k)|$ for $j=1,\ldots,m$. Note
that each $N_j(k)$ is either bounded and we can apply condition {\em (a)},
or we can apply one of the conditions {\em (c)--(e)}. The sum of all
candidates gives a series which is a uniform majorant for $\bS(n)$.\QED

\noindent We will use the following asymptotic expressions for the gamma
function. It will be convenient for us to uniformize all gamma expressions
with a linear argument in $m\to\infty$ to expressions involving only
$\Gamma(m)$. Some corollaries are formulated in less generality than
possible, for readiness of application. 
\begin{lemma} \label{gammaspm}
Let $\lambda$ be a real number, and let $\ell\in\CC$.
\begin{itemize}
\item If $\lambda>0$ then
\begin{equation} \label{gammapl}
\Gamma(\lambda\,m+\ell)\sim(2\pi)^{\frac{1-\lambda}{2}}\,\lambda^{\ell-1/2}
\;m^{\ell+\frac{\lambda-1}{2}}\;\lambda^{\lambda\,m}\;\Gamma(m)^{\lambda}
\qquad\mbox{as real $m\to\infty$}.
\end{equation}
\item If $\lambda<0$, $\ell\not\in\ZZ$, and $m$ runs through a set of real
numbers such that $\lambda\,m\in\ZZ$, then
\begin{equation} \label{gammanl}
\Gamma(\lambda\,m+\ell)\sim
\frac{(2\pi)^{\frac{1-\lambda}{2}}\,|\lambda|^{\ell-1/2}}{2\;\sin(\pi\ell)}
\;m^{\ell+\frac{\lambda-1}{2}}\;\lambda^{\lambda\,m}\;\Gamma(m)^{\lambda}
\qquad\mbox{as $m\to\infty$}.
\end{equation}
\end{itemize}
\end{lemma}
\proof The first statement follows from Stirling's asymptotic formula
\cite[Theorem 1.4.1]{specfaar}:
\begin{eqnarray*}
\frac{\Gamma(\lambda\,m+\ell)}{\Gamma(m)^{\lambda}} & \sim &
\frac{\sqrt{2\pi}\;(\lambda m+\ell)^{\lambda\,m+\ell-1/2}\;
\exp(-\lambda\,m-\ell)}{(2\pi)^{\lambda/2}\;
m^{\lambda\,m-\lambda/2}\;\exp(-\lambda\,m)}\\
&\sim&(2\pi)^{\frac{1-\lambda}{2}}\,\lambda^{\lambda\,m+\ell-1/2}
\,m^{\ell+\frac{\lambda-1}{2}}\,
\left(1\!+\!\frac{\ell}{\lambda\,m}\right)^{\ell-1/2}
\left(1\!+\!\frac{\ell}{\lambda\,m}\right)^{\lambda\,m}\,\exp(-\ell).
\end{eqnarray*}
Note that
\[
\lim_{m\to\infty}\left(1\!+\!\frac{\ell}{\lambda\,m}\right)^{\ell-1/2}=1
\qquad\mbox{and}\qquad \lim_{m\to\infty}
\left(1\!+\!\frac{\ell}{\lambda\,m}\right)^{\lambda\,m}=\exp(\ell).
\]
Formula (\ref{gammapl}) follows.

To prove the second statement we use Euler's reflection formula
\cite[Theorem 1.2.1]{specfaar}:
\begin{equation} \label{eulreflect}
\Gamma(\lambda\,m+\ell)=\frac{(-1)^{\lambda m}\;\pi}{\sin\pi\ell}
\,\frac{1}{\Gamma\left(|\lambda| m+1\!-\!\ell\right)}.
\end{equation}
Now we apply the first statement to $\Gamma\left(|\lambda|
m+1\!-\!\ell\right)$ and obtain (\ref{gammanl}). \QED

\begin{conclude} \label{gammaas}
Let $\lambda$ be a nonzero real number, and let $\ell\in\CC$. 
We assume that $m$ runs through a set of real numbers such that $\lambda\,m$
is an integer. If $\lambda<0$ then we additionally assume that
$\ell\not\in\ZZ$. Under these assumptions there is a constant $C_0\in\CC$
such that
\begin{eqnarray} \label{newgamma}
\Gamma(\lambda\,m+\ell) &\sim& C_0\;m^{\ell+\frac{\lambda-1}{2}}\;
\lambda^{\lambda\,m}\;\Gamma(m)^{\lambda}
\end{eqnarray}
as $m\to\infty$.
\end{conclude}

\begin{conclude} \label{le:gamman}
Let $\lambda$, $N$ be integers, and let $\ell\in\CC$. We assume that
$\lambda\neq 0$. If $\lambda<0$ we additionally assume that
$\ell\not\in\ZZ$. Then, as integer $m\to\infty$,
\begin{eqnarray} \label{no:gamman}
\Gamma(\lambda\,m+N+\ell) &\sim&
C_0\,(2\pi)^{\frac{1-\lambda}2}\,|\lambda|^{\ell-1/2}\;\lambda^N\;m^{N+\ell+\frac{\lambda-1}{2}}\;
\lambda^{\lambda\,m}\;\Gamma(m)^{\lambda},\qquad
\end{eqnarray}
where
\begin{equation}
C_0=\left\{\begin{array}{rl} 1, & \mbox{if } \lambda>0,\\
\frac1{2\sin\pi\ell},  & \mbox{if } \lambda<0.\end{array} \right.
\end{equation}
\end{conclude}
\proof For $\lambda<0$, the simplification is
$|\lambda|^N/\sin\pi(\ell+N)=\lambda^N/\sin\pi\ell$. \QED

\section{Other preliminary results}
\label{sec:otherprel}

Here we continue with more asymptotic formulas for the gamma function and
Pochhammer symbols. Lemma \ref{seqlemma} is used only in the auxiliary
Section \ref{gtsection}.

We introduce the following function:
\begin{equation}
\Theta(x)=\frac{1+x}{x}\,\log(1\!+\!x)-1.
\end{equation}
\begin{lemma} \label{lem:gammun}
Let $\omega(m)$ denote a real-valued function defined for large enough
$m\in\RR$, such that $\omega(m)=o(m)$ as $m\to\infty$. Then
\begin{equation} \label{lem:gammunf}
\Gamma(m+\omega(m)) \,\sim\, m^{\omega(m)}\;\exp\left(
\omega(m)\;\Theta\!\left( \frac{\omega(m)}{m} \right) \right) \; \Gamma(m)
\qquad\mbox{as}\quad m\to\infty.
\end{equation}
\end{lemma}
\proof By Stirling's asymptotic formula:
\begin{eqnarray} \label{lem:gammunf1}
\frac{\Gamma(m+\omega(m))}{\Gamma(m)} & \sim &
\frac{(m+\omega(m))^{m+\omega(m)-1/2}}{m^{m-1/2}}\,\exp(-\omega(m)) \nonumber \\
& \sim &
m^{\omega(m)}\,\left(1+\frac{\omega(m)}{m}\right)^{m+\omega(m)-1/2}\,\exp(-\omega(m))
\nonumber\\ & \sim & m^{\omega(m)}\,\exp\left(
(m+\omega(m))\,\log\left(1+\frac{\omega(m)}{m}\right)-\omega(m) \right). 
\end{eqnarray}
The result follows. \QED

\begin{conclude} \label{gammun}
Let $\omega(m)$ denote a real-valued function defined for large enough
$m\in\RR$, such that $\omega(m)=o(m)$ as $m\to\infty$. Then
\begin{equation} \label{gammunf}
\Gamma(m+\omega(m)) \sim m^{\omega(m)}\;\exp\left( \sum_{j=1}^{\infty}
\frac{(-1)^{j+1}}{j\,(j+1)}\,\frac{\omega(m)^{j+1}}{m^j}\right)\;\Gamma(m).
\end{equation}
\end{conclude}
\proof On the interval $x\in (-1,1)$ we have
\begin{equation} \label{taylorthet}
\Theta(x)=\sum_{j=1}^{\infty} \frac{(-1)^{j+1}}{j\,(j+1)}\;x^j.
\end{equation}
\QED

\noindent As a direct consequence, we obtain the following well known
asymptotics:
\begin{equation} \label{ko:pochhlm}
(\ell)_m=\frac{\Gamma(m+\ell)}{\Gamma(m)}\sim m^{\ell},\qquad\mbox{as}\quad
m\to\infty\qquad(\ell\in\CC).
\end{equation}

\begin{lemma} \label{pocchasnl}
Let $\lambda$ be a nonzero real number, and let $\ell\in\CC$. If $\lambda<0$
we additionally assume that $\ell\not\in\ZZ$. Let $N=N(m)$ denote a function
$N:\Zplus\to\Zplus$ such that $N(m)-\lambda\,m=o(m)$ as $m\to\infty$. Let
$\omega(m)$ denote the difference $N(m)-\lambda\,m$. Then there is a
constant $C_1\in\CC$ such that
\begin{equation}
(\ell)_{N(m)} \sim C_1\,m^{\ell+\frac{\lambda-1}2}\,
\lambda^{N(m)}\,m^{\omega(m)}\,\exp\left(
\omega(m)\;\Theta\!\left(\frac{\omega(m)}{\lambda\,m}\right)\right)\,\Gamma(m)^{\lambda}
\qquad\mbox{as}\quad m\to\infty.
\end{equation}
\end{lemma}
\proof We have $N(m)=\lambda\,m+\omega(m)$ and
$(\ell)_N=\Gamma(N+\ell)/\Gamma(\ell)$. Applying Corollary \ref{gammaas},
\begin{eqnarray*}
\Gamma(\lambda\,m+\omega(m)+\ell) & = & \Gamma\left( \lambda\,\left(
m+\frac{\omega(m)}{\lambda}\right)+\ell\right)\\
& \sim & C_0\,m
^{\ell+(\lambda-1)/2}\,\lambda^{\lambda\,m+\omega(m)}\;
\Gamma\!\left(m+\frac{\omega(m)}{\lambda}\right)^{\lambda}.
\end{eqnarray*}
Then we apply Lemma \ref{lem:gammun} to the last factor.\QED

\begin{conclude} \label{pocchasnbig}
Let $\lambda$ be a nonzero integer, and let $\ell\in\CC$. If $\lambda<0$, we
assume that $\ell\not\in\ZZ$. Let $\omega(m)$ denote a function
$\omega:\Zplus\to\Zplus$ such that $\omega(m)=o(m)$ as $m\to\infty$. Then
there is a constant $C_0\in\CC$ such that
\begin{equation}
(\lambda\,m+\ell)_{\omega(m)} \sim
C_0\,\lambda^{\omega(m)}\,m^{\omega(m)}\,\exp\left(
\omega(m)\;\Theta\!\left(\frac{\omega(m)}{\lambda\,m}\right)\right)
\qquad\mbox{as}\quad m\to\infty.
\end{equation}
\end{conclude}
\proof We have
$(\lambda\,m+\ell)_{\omega(m)}=\Gamma(\lambda\,m+\omega(m)+\ell)/\Gamma(\lambda\,m+\ell)$.
Lemma \ref{pocchasnl} can be applied to the numerator and the denominator.
\QED

\begin{lemma} \label{pochhntk}
Let $\lambda$ be a nonzero integer, and let $\ell\in\CC\setminus\ZZ$. Let
$\omega(m)$ denote a function $\omega:\Zplus\to\RR$ such that
$\lambda\,\omega(m)\in\ZZ$ whenever $m\in\ZZ_+$.
\begin{itemize}
\item If $\omega(m)$ is bounded as $m$ varies over $\ZZ_+$, 
then the Pochhammer symbol $(\ell)_{\lambda\,\omega(m)}$ is bounded, and it
is bounded away from zero as well.
 \item Suppose that $\omega(m)$ approaches $+\infty$ or $-\infty$ as
$m\to\infty$. Let us denote
\begin{equation} \label{sigmak}
\varepsilon=\left\{ \begin{array}{rl}
1, & \mbox{if $\omega(m)\to +\infty$ as $m\to\infty$},\\
-1, & \mbox{if $\omega(m)\to -\infty$ as $m\to\infty$}, \end{array} \right.
\end{equation}
Then there is a constant $C_0\in\CC$ such that
\begin{equation} \label{pochhas2} 
(\ell)_{\lambda\,\omega(m)} 
\sim
C_0\,|\omega(m)|^{\ell-\frac12+\frac12{\varepsilon\lambda}}\;(\varepsilon\lambda)^{\lambda\,\omega(m)}
\;\Gamma(|\omega(m)|)^{\varepsilon\lambda}.
\end{equation}
\end{itemize}
\end{lemma}
\proof If $\omega(m)$ is bounded, $\lambda\omega(m)$ achieves finitely many
values. The same holds for the Pochhammer symbol. Since $\ell\not\in\ZZ$,
the Pochhammer symbol is never zero.

If $\omega(m)\to+\infty$, by Corollary \ref{gammaas} we have
\[
\Gamma(\lambda\,\omega(m)+\ell)\sim
\widetilde{C}_0\,\omega(m)^{\ell+\frac{\lambda-1}{2}}\;\lambda^{\lambda\,\omega(m)}\;
\Gamma(\omega(m))^{\lambda}.
\]
If $\omega(m)\to-\infty$, we apply Corollary \ref{gammaas} to $\Gamma(
-\lambda\,|\omega(m)|+\ell)$. The result is
\[
\Gamma(\lambda\,\omega(m)+\ell)\sim
\widehat{C}_0\,|\omega(m)|^{\ell-\frac{\lambda+1}{2}}\;(-\lambda)^{\lambda\,\omega(m)}\;
\Gamma(|\omega(m)|)^{-\lambda}.
\]
The Pochhammer symbol grows accordingly. \QED

\begin{lemma} \label{stconverge}
Let $Z=\sum_{j=0}^\infty v_j$ be a series, $\varrho$ be a positive real
number, and let
\[
w(j)=\frac{\log |v_j|}{j^\varrho}\qquad\mbox{for}\quad j=1,2,\ldots.
\]
\begin{itemize}
\item If $\lim_{j\to\infty}w(j)=-\infty$ or $\lim_{j\to\infty}\sup w(j)<0$,
then the series $Z$ converges absolutely. \item If $\lim_{j\to\infty}\sup
w(j)=\infty$ or $\lim_{j\to\infty}\sup w(j)>0$, then the series $Z$
diverges.
\end{itemize}
\end{lemma}
\proof The case $\varrho=1$ is equivalent to the standard convergence
criteria involving $\lim\sup |v_j|^{1/j}$; see \cite{rudin}. To prove the
first statement in general, we choose a positive real number $K$ such that
$w(j)<-K$ for large enough $j$. Then $\log|v_j|<-Kj^{\varrho}<-2\,\log j$
for large enough $j$. Therefore a tail of the series $Z$ can be majorated by
$\sum j^{-2}$, so converges absolutely.

To prove the second statement we choose a positive real $K$ such that
$w(j)>K$ for arbitrary large $j$. Then $\log|v_j|>Kj^{\varrho}$, so $v_j$ is
unbounded. Hence the series $Z$ diverges. \QED

\begin{lemma} \label{seqlemma}
Let $a,b,p$ be real numbers. Assume that $a>0$ and $0<p<1$. Consider the
sequence $v_j=(aj+b)\,p^j$, with $j=0,1,2,\ldots$. Then
\[ v_{j+1}<v_j \qquad\mbox{if}\qquad j>\frac{p}{1-p}-\frac{b}{a}.
\]
\end{lemma}
\proof A straightforward computation. \QED

\section{Notation}
\label{sec:notation}

We work with a proper hypergeometric term presented in
(\ref{eq:properterm}):
\begin{equation} \label{eq4:properterm}
u(n,k)=P(n,k)\,\frac{\xi^n\theta^k}{k!}\,(b_1)_{\alpha_1n+\beta_1k}\,\cdots\,(b_p)_{\alpha_p\,n+\beta_p\,k}.
\end{equation}
We assume the conditions \refi{I}--\refi{II} presented after formula
(\ref{eq:properterm}). Then the series $\bS(n)=\sum_{k=0}^{\infty} u(n,k)$
is well defined\footnote{Slightly more generally, one may consider proper
terms with the Pocchammer symbols $(b_j)_{\alpha_jn+\beta_jk}$ replaced by
$(b_j+\alpha_j\,n)_{\beta_j\,k}$, and insist on conditions \refi{I} and
\refi{II} only when $\beta_j>0$, $\alpha_j\le 0$ respectively $\beta_j<0$,
$\alpha_j\ge 0$. Then the series $\sum_k u(n,k)$ would be well defined and
nonterminating for large enough $n$. Our results can be easily applied to
this form as well, since we can always identify
$(b_j+\alpha_j\,n)_{\beta_jk}=(b_j)_{\alpha_jn+\beta_jk}/(b_j)_{\alpha_jn}$
for large enough $n$. Example \ref{ex:hpg1} below employs reminiscent form.

A customary alternative is to write a general proper hypergeometric term in
the form
\[
u(n,k)=P(n,k)\,\frac{(b_1)_{\alpha_1n+\beta_1k}\,\cdots\,(b_p)_{\alpha_pn+\beta_pk}}
{(d_1)_{\gamma_1n+\delta_1k}\,\cdots\,(d_q)_{\gamma_qn+\delta_qk}\,k!}\,\xi^n\theta^k,
\]
where we would require $\beta_j\ge 0$ (and perhaps $\alpha_j>0$ if
$\beta_j=0$) and  $\delta_j\ge 0$ (and perhaps $\gamma_j>0$ if
$\delta_j=0$). Then notation definitions and formulas in the proof become
slightly more complicated. The first electronic
version of this article was written for this customary form.} %
for all non-negative integers $n$, and it is nonterminating for all but
possibly finitely many non-negative integers $n$. Note that we allow
$P(n,k)=0$ for infinitely many integer pairs $(n,k)$. In particular,
$P(n,k)$ may have linear factors in $n,k$.

Now we introduce a lot of notation for the expressions we need to check in order to determine
uniformly dominated convergence for the series $\bS(n)$. Example \ref{ex:hpg1} below
should be a helpful guide\footnote{As an extra guidance to our notation,
we indicate that most conditions for uniformly dominant convergence of $\bS(n)$ come from convergence conditions of several series of the following asymptotic form:
\[ 
\sim\Gamma(k)^{D_0}\,\Gamma(n)^{D^*_0}\;z_1^k\;\zeta_0^n\;k^{A_0-1}\;n^{A^*_\infty}
\,\exp\!\left(k\,\wideha{\Phi}_\infty\!\left(\frac{k}{n}\right)\right).
\] 
The variables in (\ref{no:s0}) appear in the powers of $\Gamma(k)$ and
$\Gamma(n)$ in these asymptotic forms; the expressions in
(\ref{no:a0})--(\ref{no:a1s}), (\ref{no:termsn})--(\ref{no:termsnns}) and
(\ref{no:termstr}) appear in the powers of $k$ and $n$; the expressions in
(\ref{no:z0})--(\ref{no:gt}) appear with the exponents of $k$ and $n$; the
functions in (\ref{no:phil})--(\ref{no:phii}) appear in the additional
exponent.} %
and motivator.

For those $j$ with $\beta_j\neq 0$ we introduce
\begin{equation} \label{no:hacs}
\widehat{a}_j=b_j+\frac{\beta_j-1}2.
\end{equation}

For those $j$ with $\alpha_j\neq 0$ we set
\begin{equation} \label{no:tacs}
\widetilde{a}_j=b_j+\frac{\alpha_j-1}2.
\end{equation}

Let us also introduce the following notation:
\begin{eqnarray} \label{no:s0}
D_0=-1+\sum_{j=1}^p \beta_j,\qquad D_1=\sum_{\alpha_j\neq 0}\beta_j, &&
D^*_0=\sum_{j=1}^p \alpha_j,\qquad D^*_1=\sum_{\beta_j\neq 0} \alpha_j.
\end{eqnarray}
In the sums for $D_1$ and $D^*_1$, the summation range is understood to be
over all $j$ for which $\alpha_j\ne0$ respectively $\beta_j\ne0$. In the
rest of the paper, summation or product ranges are implied by the range of
definition of involved variables and by indicated conditions. For example,
$\sum_{\alpha_j=0}\,\widehat{a}_j$ is a summation over those $j$ for which
$\beta_j>0$ 
and $\alpha_j=0$. With this convention we define:
\begin{eqnarray} \label{no:a0}
A_0&=&\;\sum_{} \widehat{a}_j+\deg_k P(n,k),\\
A^*_{\infty}\!&=&\;\sum_{} \widetilde{a}_j+\deg_n P(n,k),\\
\label{no:aos} A^*_0&=&\sum_{\alpha_j=0} \widehat{a}_j+\deg_k Q(k),\\
A_1&=& \sum \widehat{a}_j+\sum_{\beta_j=0}
\widetilde{a}_j+\deg_{\{n,k\}} P(n,k),\\
\label{no:a1s} A^*_1&=& \sum_{\alpha_j=0} \widehat{a}_j+\sum
\widetilde{a}_j+\deg_{\{n,k\}} P(n,k).
\end{eqnarray}
In (\ref{no:aos}), we denote
\begin{equation} \label{no:qk}
Q(k):=\mbox{the leading coefficient of $P(n,k)$ with respect to $n$}.
\end{equation}
Thus $Q(k)$ is a polynomial in $k$.

Now we define the function
\begin{equation} \label{no:phi}
\varphi(p)=\max_{f: \mbox{\scriptsize\,a monomial}\atop\mbox{\scriptsize of
}P(n,k) } \left( \deg_kf+p\,\deg_nf \right).
\end{equation}
Therefore $\varphi(p)$ is the degree of the polynomial $P(n,k)$ if we give
the weight $1$ to the variable $k$ and the weight $p>0$ to the variable $n$.
We have the following properties.
\begin{lemma} \label{lem:phip}
\begin{enumerate}
\item For a function $N(k):\ZZ_+\to\RR$ such that $N(k)\sim C_0\,k^p$ as
$k\to\infty$ for some nonzero constant $C_0$ and real $p>0$, we have
$P(N,k)=O(k^{\varphi(p)})$. For a general such function $N(k)$, there
is a nonzero constant $C_1$ 
such that $P(N,k)\sim C_1\,k^{\varphi(p)}$. 
  \item The function $\varphi(p)$ is a continuous piecewise linear function on the
real interval $[0,\infty)$, monotone nondecreasing. The linear slope of
$\varphi(p)$ can only increase as $p$ increases, as well.
 \item For large enough $p$, we have $\varphi(p)=p\,\deg_nP(n,k)+\deg_kQ(k)$.
\end{enumerate}
\end{lemma}
\proof The first part is clear. In particular, for any fixed $p>0$ the
second statement holds for general $C_0$.

Let ${\bf P}$ denote the Newton polygon of $P(n,k)$, that is, the convex
hull in $\RR^2$ of all half-lines from $(U,V)$ to $(U,-\infty)$ and
$(-\infty,V)$ for each monomial $k^Un^V$ of $P(n,k)$. Let
$\{(U_i,V_i)\}_{i=1}^m$ be the sequence of the vertices of ${\bf P}$,
ordered by increasing $V_i$. Then
\begin{equation}
\varphi(p)=\left\{ \begin{array}{rl} V_1p+U_1,& \mbox{if $0\le p\le
\frac{U_1-U_2}{V_2-V_1}$},\\
V_i\,p+U_i,& \mbox{for $1<i<m$ and  $\frac{U_{i-1}-U_i}{V_{i}-V_{i-1}}\le
p\le
\frac{U_i-U_{i+1}}{V_{i+1}-V_i}$},\\
V_mp+U_m,& \mbox{if $p\ge\frac{U_{m-1}-U_m}{V_{m}-V_{m-1}}$.}
\end{array}\right.
\end{equation}
The last two claims follow.\QED

\noindent Consequently, we introduce the two functions:
\begin{eqnarray} \label{no:termsn}
\psi_0(p)&=&\sum\widehat{a}_j+p\sum_{\beta_j=0}\widetilde{a}_j
+\varphi(p),\\ \label{no:termsnns}
\psi_\infty(p)&=&
\sum_{\alpha_j=0}\widehat{a}_j
+p\sum_{}\widetilde{a}_j +\varphi(p).
\end{eqnarray}
We will consider $\psi_0(p)$ on the interval $[0,1]$, and the function
$\psi_\infty(p)$ on the interval $[1,\infty)$. We have the following
properties.
\begin{lemma} \label{lem:psip}
\begin{enumerate}
\item The real parts of $\psi_0(p)$ and $\psi_{\infty}(p)$ are continuous
piecewise linear functions on the real interval $[0,\infty)$. Their linear
slopes can only increase as $p$ increases.
 \item On any interval $[U,V]\subset [0,\infty)$, the real parts of $\psi_0(p)$
and $\psi_\infty(p)$ achieve their maximum on $[U,V]$ at an end point, $U$
or $V$.
 \item Let $(U,V)$ be a subinterval $[0,\infty)$, so possibly $V=\infty$.
If the linear slope of $\mbox{\rm Re }\psi_0(p)$ or $\mbox{\rm Re
}\psi_{\infty}(p)$ is zero or negative as $p\to V$ from the left, then the
supremum of $\mbox{\rm Re }\psi_0(p)$ or $\mbox{\rm Re }\psi_{\infty}(p)$ on
$(U,V)$ is approached as $p\to U$.
 \item $\psi_0(0)=A_0$, $\psi_0(1)=A_1$, and $\psi_\infty(1)=A^*_1$.
 \item For large enough $p$, we have $\psi_\infty(p)=p\,A^*_\infty+A_0^*$.
\end{enumerate}
\end{lemma}
\proof The first part follows from Lemma \ref{lem:phip} {\em (ii)}. Since
the slopes can only increase, on each interval $[U,V]$ the real parts of
$\psi_0(p)$ and $\psi_{\infty}(p)$ are either monotone functions, or there
is one locally extremal value inside the interval and that value is a local
minimum. This shows the second part. In part {\em (iii)}, the function
$\mbox{\rm Re }\psi_0(p)$ or $\mbox{\rm Re }\psi_{\infty}(p)$ does not
increase on $(U,V)$. The last two parts are straightforward. \QED

\noindent Let us define the set
\begin{equation} \label{omegadef}
\Omega=\left\{-\frac{\beta_j}{\alpha_j}\right\}_{\alpha_j\neq 0},
\end{equation}
and the family of polynomials
\begin{equation}
P^{\star}_\lambda(n,k):=P(\lambda k+n,k).
\end{equation}
We assume that the polynomial $P^{\star}_\lambda(n,k)$ is expanded whenever
we implicitly use it for some $\lambda$. Similarly as in (\ref{no:phi}), we
define the family of functions
\begin{equation} \label{no:phist}
\varphi^{\star}_{\lambda}(p)=\max_{f:\mbox{\scriptsize\,a
monomial}\atop\mbox{\scriptsize of }P^{\star}_{\lambda}(n,k) } \left(
\deg_kf+p\,\deg_nf \right).
\end{equation}
We may need to consider these functions on the interval $p\in[0,1]$.
\begin{lemma} \label{lem:phipst}
\begin{enumerate}
\item We have $\varphi^{\star}_{\lambda}(1)=\varphi(1)$, and
$\varphi^{\star}_{\lambda}(p)\le\varphi(1)$ for $p\in[0,1]$. If
\begin{equation} \label{ine:pkk}
\deg_k P(\lambda\,k,k)=\deg_{\{n,k\}} P(n,k),
\end{equation}
then $\varphi^{\star}_{\lambda}(p)=\varphi(1)$ for any $p\in[0,1]$.
 \item The function $\varphi^{\star}_{\lambda}(p)$ is
a continuous piecewise linear function on the real interval $[0,1]$,
monotone nondecreasing. The linear slope of $\varphi^{\star}_{\lambda}(p)$
can only increase as $p$ increases, as well.
 \item Let $N(k):\ZZ_+\to\RR$ denote a function such that
$N(k)\sim\lambda\,k+C_0\,k^p$ as $k\to\infty$ for some nonzero constants
$C_0\neq 0$, $\lambda\neq 0$ and $p\in[0,1)$. Then $P(N,k)\sim C_1
k^{\varphi^{\star}_{\lambda}(p)}$ for some nonzero constant $C_1$.
\end{enumerate}
\end{lemma}
\proof For the first part, note that $\deg_{\{n,k\}}
P^{\star}_{\lambda}(\lambda\,k,k)=\deg_{\{n,k\}} P(n,k)$. If (\ref{ine:pkk})
is satisfied, then the coefficient of $P^{\star}_\lambda(n,k)$ to
$k^{\varphi(1)}$ is nonzero.
 (Non-generically, we may have $\deg_k P(\lambda\,k,k)<\deg_{\{n,k\}} P(n,k)$.)

The other two parts follow similarly as parts {\em (ii)} and {\em (i)} of
Lemma \ref{lem:phip}, respectively.   \QED

\noindent We introduce a variation of $\psi_0(p)$ as well:
\begin{eqnarray} \label{no:termstr}
\psi^{\star}_{\lambda}(p)=\sum_{\alpha_j\lambda+\beta_j\neq0}\!
\left(\widehat{a}_j+\,\frac{\alpha_j}2\right)+\sum_{\beta_j=0}
\widetilde{a}_j\,
+p\!\sum_{\alpha_j\lambda+\beta_j=0} \widehat{a}_j
+\,\varphi^{\star}_{\lambda}(p).
\end{eqnarray}
Note that the linear coefficient to $p$ is zero if $\lambda\not\in\Omega$.
\begin{lemma} \label{lem:psipstr}
\begin{enumerate}
\item For generic $\lambda$, the function $\psi^{\star}_{\lambda}(p)$ is a
constant:
\begin{equation} \label{no:b1str}
\psi^{\star}_{\lambda}(p)=A_1+\frac{D^*_1}2 \;=\;A_1^*+\frac{D_1}2.
\end{equation}
\item For any $\lambda$, the real part of $\psi^{\star}_{\lambda}(p)$ is a
continuous piecewise linear function on the real interval $[0,1]$. Its
linear slope can only increase as $p$ increases.
 \item On any interval $[U,V]\subset [0,1]$, the real part of
$\psi^{\star}_{\lambda}(p)$ achieve its maximum on $[U,V]$ at an end point,
$U$ or $V$. If the linear slope of $\mbox{\rm Re }\psi^{\star}_{\lambda}(p)$
is zero or negative as $p\to V$ from the left, then the supremum of
$\mbox{\rm Re }\psi^{\star}_{\lambda}(p)$ on $(U,V)$ is approached as $p\to
U$.
\end{enumerate}
\end{lemma}
\proof In the first part, the generic $\lambda$ are those
$\lambda\not\in\Omega$ 
which satisfy (\ref{ine:pkk}). Other two parts follow similarly as parts
{\em (i)--(iii)} of Lemma \ref{lem:psip}.\QED

\begin{example} \label{ex:hpg1} \rm Consider the hypergeometric series
\begin{equation} \label{eq:canpropterm}
\hpg{s}{r}{\wideha{b}_1+\alpha_1n,\,\ldots,\,\wideha{b}_S+\alpha_Sn,\,
\wideha{b}_{S+1},\ldots,\wideha{b}_s}
{\wideha{d}_1+\gamma_1n,\,\ldots,\,\wideha{d}_R+\gamma_Rn,\,
\wideha{d}_{R+1},\ldots,\wideha{d}_r}{\,Z},
\end{equation}
where $\alpha_i$'s and $\gamma_i$'s are nonzero integers. To put the
hypergeometric series in the form (\ref{eq4:properterm}), we may rewrite it
as 
\[ 
\sum_{k=0}^{\infty}
\prod_{j=1}^S\frac{(\wideha{b}_j)_{\alpha_jn+k}}
{(\wideha{b}_j)_{\alpha_jn}}
\prod_{j=1}^R\frac{(\wideha{d}_j)_{\gamma_jn}}
{(\wideha{d}_j)_{\gamma_jn+k}} 
\,\frac{\prod_{j=S+1}^s(\wideha{b}_j)_{k}}
{\prod_{j=R+1}^r(\wideha{d}_j)_{k}}\,\frac{Z^k}{k!},
\] 
and then move the Pochhammer symbols in the denominators by using
$$
\frac1{(\wideha{b}_j)_{\alpha_jn}}=(1-\wideha{b}_j)_{-\alpha_jn}(-1)^{\alpha_jn},\quad
\frac1{(\wideha{d}_j)_{\gamma_jn+k}}=(1-\wideha{d}_j)_{-\gamma_jn-k}(-1)^{\gamma_jn+k},
\quad  \mbox{etc.}
$$

In the setting of (\ref{no:hacs}) and (\ref{no:tacs}), the set of all $\widehat{a}_j$'s is 
$$\{b_j\}_{j=1}^s\cup\{-d_j\}_{j=1}^r,$$ and the set of all $\widetilde{a}_j$'s is 
$$\textstyle\left\{b_j+\frac{\alpha_j-1}2\right\}_{j=1}^S\cup
\left\{-b_j+\frac{1-\alpha_j}2\right\}_{j=1}^S\cup\left\{d_j+\frac{\gamma_j-1}2\right\}_{j=1}^R
\cup\left\{-d_j+\frac{1-\gamma_j}2\right\}_{j=1}^R.$$
We have:
\begin{eqnarray}
D_0=s-r-1,\qquad D^*_0=0,\qquad D_1=S-R,  &&
D^*_1=\sum_{j=1}^S\alpha_j-\sum_{j=1}^R\gamma_j,\\
A^*_0=A^*_1=\sum_{j=S+1}^s b_j-\sum_{j=R+1}^r d_j,&&
A_0=\sum_{j=1}^s b_j-\sum_{j=1}^r d_j, \\
A_1=A^*_1+\frac{D_1-D_1^*}2,\hspace{18pt} && A^*_\infty=0.
\end{eqnarray}
The functions $\varphi(p)$ and $\varphi^\star_\lambda(p)$ in 
(\ref{no:phi}) and (\ref{no:phist}) are identically equal to 0, and
$\psi_\infty(p)= A^*_1$. We have
\begin{eqnarray}
\psi_0(p)&\!\!=\!\!& A^*_1+
(1-p)\left(\sum_{j=1}^S b_j-\sum_{j=1}^R d_j\right)+
\frac{p}2\left(D_1-D_1^*\right),\\
\psi^\star_\lambda(p)&\!\!=\!\!& A^*_1+\frac{D_1}2-
(1-p)\!\left(\sum_{\alpha_j\lambda=-1}\!\!b_j
-\sum_{\gamma_j\lambda=-1}\!\!d_j\right)
+\frac{D_{\lambda}}{2\,\lambda}.
\end{eqnarray}
where $D_\lambda$ is the difference between the number of $\alpha_j$'s equal to $-1/\lambda$
and the number of $\gamma_j$'s equal to $-1/\lambda$.
This example is continued below as Example \ref{ex:hpg2}.
\end{example}

\section{Further notation}
\label{sec:enotation}

The notation of the previous Section adds up the the parameters
$\alpha_j,\beta_j,b_j$. Here we introduce some ``multiplicative" notation.
Recall the convention $0^0=1$.

We introduce the following constants:
\begin{eqnarray} \label{no:z0}
z_0=\theta\;{\prod \beta_j^{\beta_j}},
 & \displaystyle\quad z_1=\theta\,{\prod_{\alpha_j\neq 0}\alpha_j^{\beta_j}}\,
 {\prod_{\alpha_j=0}\beta_j^{\beta_j}},\quad
 & z_\infty= \theta\,\prod_{\alpha_j\neq 0}\alpha_j^{\beta_j},\\
\zeta_0=\xi\;{\prod \alpha_j^{\alpha_j}},& \displaystyle \zeta_1=\xi\,
{\prod_{\beta_j\neq0}\beta_j^{\alpha_j}}\,{\prod_{\beta_j=0}\alpha_j^{\alpha_j}}.
\end{eqnarray}
Besides, we define the function
\begin{equation} \label{no:gt}
g(t)=|\theta|\,|\xi|^t\;{\prod|\beta_j+\alpha_jt|^{\beta_j+\alpha_jt}}.
\end{equation}
We have the following properties of $g(t)$.
\begin{lemma} \label{gtgenlemma}
\begin{enumerate}
\item The function $g(t)$ is continuous on the whole real axis. It can be
expressed as follows:
\begin{equation}  \label{eq:gt}
g(t)=|z_1|\,|\zeta_0|^t\,
{\prod_{\alpha_j\neq0}\left|\,t+\frac{\beta_j}{\alpha_j}\,\right|^{\beta_j+\alpha_jt}}.
\end{equation}
 \item $g(t)$ is continuously differentiable on $\RR\setminus\left(\{0\}
\cup\Omega\right)$, and
\begin{eqnarray} \label{gtderiv}
\exp\frac{g'(t)}{g(t)}&=&|\xi|\;\exp(D^*_0)\; 
{\prod|\beta_j+\alpha_jt|^{\alpha_j}}\\
\label{no:gtderiv} 
&=&|\zeta_0|\,\exp(D^*_0)\;
{\prod_{\alpha_j\neq0}\left|\,t+\frac{\beta_j}{\alpha_j}\,\right|^{\alpha_j}}.
\end{eqnarray}
 \item A point $\lambda\in\{0\}\cup\Omega$ is a genuine point of discontinuity
of the derivative $g'(t)$ if and only if
$\sum_{\beta_j+\alpha_j\lambda=0}\alpha_j\neq0$. If this is the case, then
the tangent line to $g(t)$ approaches the vertical line as $t\to\lambda$.
 \item $g(0)=|z_0|$.
 \item $g(t)\sim|z_1|\,\exp(D_1)\;|\zeta_0|^t\;
t^{\,D^*_0\,t+D_1}$ as $t\to\infty$.
\end{enumerate}
\end{lemma}
\proof Consider the function
\begin{equation}
f(x)=\left\{
\begin{array}{rl} |x|^x, & \mbox{if } x\neq 0,\\
1, & \mbox{if } x=0. \end{array}\right.
\end{equation}
We can write $f(x)=\exp(x\,\log|x|)$ for nonzero $x$. It is a standard
analysis exercise that $f(x)$ is a continuous function. Since
$f'(x)\!=\!(1\!+\log|x|)\,f(x)$, the function $f(x)$ is continuously
differentiable on $\RR\setminus\{0\}$. Expressions
(\ref{eq:gt})--(\ref{no:gtderiv}) routinely follow.

For part {\em(iii)}, we compute that as $t\to\lambda$,
\begin{eqnarray} \label{as:dgt}
\exp\frac{g'(t)}{g(t)}\sim |\xi|\,
{\prod_{\beta_j+\alpha_j\lambda\neq0}\!|\beta_j+\alpha_j\lambda|^{\alpha_j}}\,
{\prod_{\beta_j+\alpha_j\lambda=0}\!|\alpha_j|^{\alpha_j}}\,\exp(D^*_0)\,
|t-\lambda|^{\sum_{\beta_j+\alpha_j\lambda=0}\,\alpha_j}.
\end{eqnarray}
Hence, as $t\to\lambda$,
\begin{equation}
g'(t)\sim\left( C_0+
\log|t-\lambda|\sum_{\beta_j+\alpha_j\,\lambda=0}\alpha_j
 \right)g(\lambda)
\end{equation}
for a constant $C_0$. 
Part {\em(iii)} is evident.

Part {\em(iv)} is obvious. To show the asymptotic expression of part
{\em(v)}, we use (\ref{eq:gt}) to derive
\begin{equation}
g(t) 
=|z_1|\;|\zeta_0|^t\;t^{D^*_0\,t+D_1}\; \prod_{\alpha_j\neq 0}
\left|1+\frac{\beta_j}{\alpha_jt}\right|^{\alpha_jt}.
\end{equation}
Whether $\lambda>0$ or $\lambda<0$, we have
$\left(1+\frac{\ell}{\lambda\,t}\right)^{\lambda\,t}\to\exp(\ell)$ as
$t\to\infty$.\QED

\noindent For completeness, one can compute that
\begin{eqnarray} \label{eq:ldgt0}
\exp\frac{g'(t)}{g(t)}&\sim&|\zeta_1|\,\exp(D^*_0)\; t^{D^*_0-D^*_1} \qquad
\mbox{as } t\to+0.\\
\exp\frac{g'(t)}{g(t)}&\sim&|\zeta_0|\,\exp(D^*_0)\; t^{D^*_0} \qquad
\mbox{as } t\to\infty.
\end{eqnarray}
The first expression is a special case of (\ref{as:dgt}). The function $g(t)$ is
examined more closely in Section \ref{gtsection}.

At the last, we introduce the family of functions
\begin{eqnarray} \label{no:phil}
\Phi_{\lambda}(x)=\sum_{\alpha_j\lambda+\beta_j\neq 0}
\frac{\alpha_j^2\,x}{\alpha_jx+\alpha_j\lambda+\beta_j}.
\end{eqnarray}
In particular,
\begin{eqnarray} \label{no:phi0}
\Phi_0(x)=\sum_{\beta_j\neq 0} \frac{\alpha_j^2\,x}{\alpha_jx+\beta_j}.
\end{eqnarray}
We also introduce
\begin{eqnarray} \label{no:phii}
\Phi_\infty(x)=\sum_{\alpha_j\neq 0} \frac{\beta_j^2\,x}{\beta_jx+\alpha_j}.
\end{eqnarray}
This is almost all notation we will need to describe the constants we have
to check to determine uniformly dominated convergence of $\bS(n)$.

\begin{example} \label{ex:hpg2}  \rm Continuing Example \ref{ex:hpg1}, we have:
\begin{eqnarray}
z_0=Z,\qquad  z_1=z_\infty=Z\,\frac{\prod_{j=1}^S\alpha_j}{\prod_{j=1}^R\gamma_j},\qquad
\zeta_0=1, && \zeta_1=\frac{\prod_{j=1}^R\gamma_j^{\gamma_j}}
{\prod_{j=1}^S\alpha_j^{\alpha_j}},\\
g(t)=|Z|\,\prod_{j=1}^S \frac{|\alpha_jt+1|^{\alpha_jt+1}}{|\alpha_jt|^{\alpha_jt}}
\prod_{j=1}^R \frac{|\gamma_jt|^{\gamma_jt}}{|\gamma_jt+1|^{\gamma_jt+1}}.
\end{eqnarray}
We also have 
\begin{eqnarray} 
\Phi_\lambda(x)&=&\sum_{\alpha_j\lambda\neq-1}
\frac{\alpha_j^2\,x}{\alpha_jx+\alpha_j\lambda+1}-
\sum_{\gamma_j\lambda\neq-1}
\frac{\gamma_j^2\,x}{\gamma_jx+\gamma_j\lambda+1}-
\frac{D_1^*\,x}{x+\lambda},\nonumber\\
\Phi_0(x)&=&\sum_{j=1}^S\frac{\alpha_j^2\,x}{\alpha_jx+1}-
\sum_{j=1}^R\frac{\gamma_j^2\,x}{\gamma_jx+1},\\
\Phi_\infty(x)&=&\sum_{j=1}^S\frac{x}{x+\alpha_j}-
\sum_{j=1}^R\frac{x}{x+\gamma_j}.\nonumber
\end{eqnarray} 
\end{example}

Somewhat more generally, if we multiply the hypergeometric series in (\ref{eq:canpropterm})
by the gamma factor
\begin{equation} \label{eq:gammaff}
\frac{\Gamma(n+a_1)\,\cdots\,\Gamma(n+a_K)}
{\Gamma(n+c_1)\,\cdots\,\Gamma(n+c_L)},
\end{equation}
then the following values and functions in Examples \ref{ex:hpg1} and \ref{ex:hpg2} change:
$A_\infty^*$, $A_1$, $A_1^*$ and $\psi^\star_{\lambda}(p)$ are increased by $\sum a_i-\sum c_i$;
the functions $\psi_0(p)$ and $\psi_{\infty}(p)$ are increased by $p\left(\sum a_i-\sum c_i\right)$;
we get $D_0^*=K-L$; the function $\Psi_{\lambda}(x)$ is increased by $(K-L)x/(x+\lambda)$; 
and the function $g(t)$ gets multiplied by $t^{(K-L)\,t}$. In particular, if $K=L$ and 
$\sum a_i=\sum c_i$, then none of the introduced values and functions changes.

\section{The main result}
\label{sec:theresult}

Our main result is the following.
\begin{theorem}  \label{th:uctheorem}
The series $\bS(n)=\sum_k u(n,k)$ defined by $(\ref{eq4:properterm})$ is
uniformly bounded by an absolutely convergent series only if the following
restrictions are satisfied:
\begin{enumerate}
\item $D_0\le 0$ and $D^*_0\le 0$.
 \item If $D_0=0$ then one of the following two conditions must hold:
\begin{itemize}
\item $|z_0|<1$.
 \item $|z_0|=1$, $\re{A_0}<0$ and $D^*_1\le0$.
\end{itemize}
\item If $D^*_0=0$ then one of the following three conditions must hold:
\begin{itemize}
\item $|\zeta_0|<1$.
 \item $|\zeta_0|=1$, $\re{A_\infty^*}<0$ and $D_1\le 0$.
 \item $\zeta_0=1$, $A_\infty^*=0$ and $D_1\le 0$.
\end{itemize}
\end{enumerate}
These conditions are sufficient for uniformly dominated convergence if
$D_0<0$ or $D^*_0<0$. Otherwise, that is when 
\begin{equation}
D_0=0\qquad \mbox{and}\qquad D^*_0=0,
\end{equation}
the series $\bS(n)$ are bounded by an absolutely convergent series if and
only if:
\begin{enumerate}
\item[(iv)] $g(t)\le 1$ for all $t>0$.
 \item[(v)] For those $t>0$ which
satisfy $g(t)=1$, we have $\re{\psi^{\star}_{t}(0)}<0$,
\begin{equation} \label{rpartv}
\sum_{\alpha_jt+\beta_j=0}\alpha_j= 0,
\end{equation}
and one of the following two conditions holds:
\begin{itemize}
\item $\Phi_{t}(x)\equiv 0$; and $\re{\psi^{\star}_{t}(1)}\le 0$.
 \item $\Phi_{t}(x)=v_mx^m+O(x^{m+1})$ around $x=0$, where $m$ is a
positive odd 
integer, $v_m<0$; and $\re{\psi^{\star}_{t}\!\left(\frac{m}{m+1}\right)}<0$.
\end{itemize}
\item[(vi)] If $|z_0|=1$ and $D^*_1=0$, then one of the following conditions
holds:
\begin{itemize}
\item $|\zeta_1|<1$.
 \item $|\zeta_1|=1$; $\Phi_{0}(x)\equiv 0$; and $\re{A_1}\le 0$.
\item $|\zeta_1|=1$; $\Phi_{0}(x)=v_mx^m+O(x^{m+1})$ around $x=0$ for some
positive integer $m$ and negative real $v_m$; and
$\re{\psi_0\!\left(\frac{m}{m+1}\right)}<0$.
\end{itemize}
\item[(vii)] 
If $|\zeta_0|=1$ and $D_1=0$,  then one of the following conditions holds:
\begin{itemize}
\item $|z_1|<1$.
 \item $|z_1|=1$; $\Phi_{\infty}(x)\equiv 0$; and $\re{A^*_1}<0$.
 \item $|z_1|=1$; $\Phi_{\infty}(x)\equiv 0$; $\re{A^*_1}=0$; and either
$\re{A^*_\infty}<0$ or
\begin{equation} \label{nv:pnk}
\deg_{\{n,k\}}P(n,k)>\deg_nP(n,k)+\deg_kQ(k).
\end{equation}
 \item $|z_1|=1$; $\Phi_{\infty}(x)=v_mx^m+O(x^{m+1})$ around $x=0$ for some
positive integer $m$ and negative real $v_m$;
$\re{\psi_{\infty}\!\left(\frac{m+1}{m}\right)}<0$; and if ${A^*_\infty}=0$
then $\re{A^*_0}<0$.
\end{itemize}
\end{enumerate}
If these conditions are satisfied, then the limit series
$\lim_{n\to\infty}\bS(n)$ is equal to:
\begin{itemize}
\item If $D^*_0<0$, $|\zeta_0|<1$ or $\re{A_\infty^*}<0$, then $0$.%
\item If $D^*_0=0$, $\zeta_0=1$, $A_\infty^*=0$ and $D_1<0$, then $H_0\,Q(0)$.
 \item If $D^*_0=0$, $\zeta_0=1$, $A_\infty^*=0$ and $D_1=0$, then
\begin{equation} \label{glimzn}
\lim_{n\to\infty}\bS(n)=H_0\,\sum_{k=0}^{\infty}
Q(k)\,\frac{z_\infty^{\,k}}{k!}\, {\prod_{\alpha_j=0}\,(b_j)_{\beta_j\,k}},
\end{equation}
\end{itemize}
where $H_0$ is the following constant:
\begin{equation} \label{no:h0}
H_0=\left(2\pi\right)^{\sum_{\alpha_j\neq0}\frac{1-\alpha_j}2}
\prod_{\alpha_j\neq0}\frac{|\alpha_j|^{b_j-\frac12}}{\Gamma(b_j)}\left/
{\prod_{\alpha_j<0} 2\sin\pi b_j}\right..
\end{equation}
\end{theorem}

\noindent We prove this Theorem in the following Section.

Here we make a few comments. We reformulate some conditions, or indicate
some possible or typically effective simplifications. We keep some
redundancy in notation or formulation to make the proof more smooth, or to
make nontypical complications better understandable. 
\renewcommand{\labelenumi}{\it(\Roman{enumi})}
\begin{enumerate}
\item Uniformly dominated convergence of ${\cal U}(n)$ does not change if it is multiplied
by the gamma term in (\ref{eq:gammaff}) with $K=L$ and $\sum a_i=\sum c_i$, because
that does not change any of the introduced values and functions.
\item Condition (\ref{nv:pnk}) means that among the monomials of $P(n,k)$ of
the highest degree in $n,k$ there are no monomials of the highest degree in
$n$. Recall that $Q(k)$ is defined in (\ref{no:qk}).
\item Equality (\ref{rpartv}) is trivially satisfied if $t\not\in\Omega$.
Recall that the set $\Omega$ is defined in (\ref{omegadef}).
 \item Let $B_1$ denote the constant 
in (\ref{no:b1str}). If $t\not\in\Omega$, then the function
$\psi^{\star}_t(p)$ in condition {\em (v)} is rather simple:
\begin{equation}
\psi^{\star}_t(p)=B_1-\varphi(1)+\varphi^{\star}_t(p).
\end{equation}
From part {\em (i)} of Lemma \ref{lem:phipst} it follows that
$\psi^{\star}_t(p)\le B_1$ for $p\in[0,1]$, and $\psi^{\star}_t(1)=B_1$. For
generic $p$ we have $\psi^{\star}_t(p)=B_1$ for $p\in[0,1]$.

If $t\not\in\Omega$, $g(t)=1$ and $\Phi_t(x)\equiv0$, the condition {\em(v)}
can be replaced by the following restriction: either $\re{B_1}<0$, or
$\re{B_1}=0$ and \mbox{$\deg_k P(tk,k)<\deg_{\{n,k\}}P(n,k)$}. If we can
apply this simpler restriction to at least one $t\not\in\Omega$, then the
points $t\not\in\Omega$ with $g(t)=1$, $\Phi_t(x)\not\equiv0$ can only
strengthen $\re{B_1}=0$ to $\re{B_1}<0$ and add conditions on $\Phi_t(x)$.
 \item Condition (\ref{rpartv}) is equivalent to $\sum_{\alpha_jt+\beta_j=0}\beta_j=0$.
For the purposes of Theorem \ref{th:uctheorem}, we can replace the third summation
in definition (\ref{no:termstr}) of $\psi^{\star}_t(p)$ by the same
conditional summation of $\widetilde{a}_j$'s or $b_j-\frac12$, because the function
$\psi^{\star}_t(p)$ is relevant only if (\ref{rpartv}) holds.
 \item All conditions of Theorem \ref{th:uctheorem} can be checked
algorithmically, if all parameters in (\ref{eq:properterm}) are given numerically.
The only less straightforward part is checking the condition $g(t)\le 1$ for $t>0$,
and identifying the points with $g(t)=1$.
We consider this issue in Section \ref{gtsection}.
 \item Suppose that the polynomial $P(n,k)$ has a linear factor
$\widetilde{\alpha}n+\widetilde{\beta}k+\ell$ with
$\widetilde{\alpha},\widetilde{\beta}\in\ZZ$. The linear factor can be
expressed as $\ell\,(\ell+1)_{\widetilde{\alpha}n+\widetilde{\beta}k}\big/
(\ell)_{\widetilde{\alpha}n+\widetilde{\beta}k}$. Notice that all
conditions, in particular {\em (v)}, are stable if we rewrite expression
(\ref{eq4:properterm}) of $u(n,k)$ by replacing in $P(n,k)$ the linear
factor by the constant $\ell$, and appending the two Pochhammer terms to the
product of Pochhammer symbols.
 \item 
The polynomial $Q(k)$ occurs only in (\ref{nv:pnk}) and in the expressions
for $\lim_{n\to\infty}\bS(n)$. The constant $z_{\infty}$ occurs only in
(\ref{glimzn}). The constant $A_0^*$ occurs in the last case of condition
{\em (vii)}. The set $\Omega$ does not explicitly appear in the formulated
Theorem.
 \item Notice that $H_0=\prod_{\alpha_j<0}\Gamma(1-b_j)\big/
\prod_{\alpha_j>0}\Gamma(b_j)$ if all nonzero $\alpha_j$'s are equal to 1 or
$-1$.
\end{enumerate}
When discussing the function $g(t)$ in Section \ref{gtsection}, we add a few more observations
\refi{X}--\refi{XIV} to this list. 

\section{Examples}
\label{examplez}

\begin{example} \label{ex:gauss} \rm
Consider the hypergeometric series $S_1(n)=\hpg{2}{1}{a+\alpha\,n,\,b}{c+\gamma\,n}{1}$,
with $\alpha,\gamma$ nonzero integers. 
Following Examples \ref{ex:hpg1} and \ref{ex:hpg2}, we have:
\begin{eqnarray}
D_0=D_1=D^*_0=0,\qquad D^*_1=\alpha-\gamma,&& z_0=\zeta_0=1,\qquad
z_1=\frac{\alpha}{\gamma},\\
A_0=a+b-c,\qquad A^*_\infty=0.
\end{eqnarray}
Of the conditions \refi{i}--\refi{iii}, we have \refi{i} and the third option of \refi{iii} satisfied.
The second option of \refi{ii} requires $\re{A_0}<0$ and $D^*_1\le 0$ for uniformly dominant
convergence of $S_1(n)$. We notice right away that condition \refi{vii} requires $|z_1|\le 1$. 
This gives the following inequalities: 
\begin{equation}
\re{(c-a-b)}<0,\qquad \alpha\le\gamma,\qquad |\alpha|\le|\gamma|.
\end{equation}
The last two inequalities mean that either $|\alpha|\le\gamma$, or $\alpha=\gamma<0$.
We consider the cases $|\alpha|<\gamma$, $\alpha=-\gamma$, $\alpha=\gamma>0$ and
$\alpha=\gamma<0$ separately. For conditions \refi{iv} and \refi{v}, we have:
\begin{equation} \label{eq:ex71gt}
g(t)=\frac{|1+\alpha\,t|^{1+\alpha\,t}\,|\gamma\,t|^{\gamma\,t}}
{|\alpha\,t|^{\alpha\,t}\,|1+\gamma\,t|^{1+\gamma\,t}}.
\end{equation}

If $|\alpha|<\gamma$, then it is straightforward to check that $g(t)$ monotonically 
decreases from 1 to $|\alpha/\gamma|$ as $t$ varies over $(0,\infty)$;
see also part \refi{iii} of Lemma \ref{gtsplemma} below. Then Theorem \ref{th:uctheorem}
gives no additional conditions, since $D^*_1<0$ and $|z_1|<1$. The series $S_1(n)$
converges uniformly, and the limit series has the value $(1-\alpha/\gamma)^{-b}$.

If $\alpha=-\gamma$ (and $\gamma>0$), then $g(t)<1$ for $t\in (0,\infty)$ as well;
see part \refi{iv} of Lemma \ref{gtsplemma}. It remains to check condition \refi{vii}; 
only the last option can be satisfied, since $\Psi_\infty(x)\not\equiv 0$. We compute
$\psi_\infty(p)=A_0^*=b$, and
\begin{equation}
\Psi_\infty(x)=\frac{x}{x-\gamma}-\frac{x}{x+\gamma}=-\frac2{\gamma}\,x+O(x^3).
\end{equation}
The restriction on $\Psi_\infty(x)$ is satisfied. But we have the additional condition $\re{b}<0$.
The limit series has the value $2^{-b}$.

If $\alpha=\gamma>0$, then $g(t)=1$. Then $\xi_1=1$, $A_1=A^*_1=b$ and 
\begin{eqnarray}
\Psi_\lambda(x)\equiv0 &&
\mbox{for $\lambda=0$, $\lambda=\infty$ or $\lambda\in(0,\infty)$},\\
\psi^\star_\lambda(p)=b  &&
\mbox{for $\lambda\in(0,\infty)$}.
\end{eqnarray}
Besides, $P(n,k)=Q(k)=1$. The first option of \refi{v} and the second options of \refi{vi}, \refi{vii}
are relevant. We get the additional condition $\re{b}<0$ here as well. The limit series
has the value $0$.

If $\alpha=\gamma<0$, a crucial difference from the previous case is that the set $\Omega$ is
nonempty. In particular, for $t=-1/\gamma$ we have $\psi^\star_t(p)=b+(1-p)(c-a)$.
Condition \refi{v} applied to this $t$ immediately gives $\re{b}<\re{(a-c)}$. This is an additional condition
to $\re{b}<0$ and $\re{(c-a-b)}>0$.

If we assume $\alpha=0$ and $\gamma>0$ for the same series  $S_1(n)$,
then the series converges uniformly if $\re{(c-a-b)}>0$, because
\begin{eqnarray}
D_0=D^*_0=0, \quad D_1=-1,\quad D^*_1=-\gamma,&& z_0=\zeta_0=1,\\
A_0=a+b-c,\qquad A^*_\infty=0,
\end{eqnarray}
and $g(t)<1$ for all $t>0$. See also \cite[Section 2]{koornwinder}. On the other hand,
if $\alpha=0$ and $\gamma<0$, then $D_1^*>0$ in condition \refi{ii}. Then the limit series
do not converge.

In conclusion, the hypergeometric series 
$\hpg{2}{1}{a+\alpha\,n,\,b}{c+\gamma\,n}{1}$ converges uniformly
in the following situations:
\begin{itemize}
\item $|\alpha|<\gamma$ and $\re{(c-a-b)}>0$. 
\item $|\alpha|=\gamma$, and 
$\re{(c-a-b)}>0$ and $\re{b}<0$.
\item $\alpha=\gamma<0$, and $\re b<-\left|\re{(c-a)}\right|$.
\end{itemize}
Notice that the series $\hpg{2}{1}{a+\alpha\,n,\,b+\beta\,n}{c+\gamma\,n}{1}$ with 
$\alpha\beta\neq0$ cannot converge uniformly, because the termwise limit does not exist.
In the setting of Theorem \ref{th:uctheorem}, we would have $D_1>0$, so condition \refi{iii}
would not be satisfied.
\end{example}

\begin{example} \rm
Consider the hypergeometric series $S_2(n)=\hpg{2}{1}{a+\alpha\,n,\,b}{c+\gamma\,n}{-1}$. Compared with the previous example, only the values $z_0$, $z_1$ (and $z_{\infty}$) change --- they get multiplied by $-1$. That does not change any of the conditions of Theorem \ref{th:uctheorem}. Hence we have the same uniform convergence cases as in the previous  example. 

In particular, the well poised series $\hpg{2}{1}{a+2n,\,b}{1+a-b+2n}{-1}$ converges uniformly if $\re{b}<0$, as used in \cite{gauthier}. The limit series has the value $2^{-b}$. But the series $\hpg{2}{1}{a-2n,\,b}{1+a-b-2n}{-1}$ fails the test of Theorem \ref{th:uctheorem}, because of a contradictory condition
$\re{b}<\re{(b-1)}$. 

Similarly, the well poised series $\hpg{2}{1}{a,\,b-n}{1+a-b+n}{-1}$ converges uniformly if 
$\re{b}<\frac12$ and $\re{a}<0$. The limit series has the value $0$. Recall that well poised
$\hpgo{2}{1}(-1)$ series 
can be evaluated using Kummer's formula \cite[Cor. 3.1.2]{specfaar}.
\end{example}

\begin{example} \rm
Consider the hypergeometric series $S_3(n)=\hpg{2}{1}{a+\alpha\,n,\,b,\,c}{d+\gamma\,n,\,f}{1}$. 
Compared with Example \ref{ex:gauss}, we have $A_0=a+b+c-d-f$ and
$\psi_\infty(p)=A_0^*=A_1=A_1^*=b+c-f$. The functions $\psi_\infty(p)$ and $\psi^\star_\lambda(p)$
are different similarly. With the same reasoning as in Example \ref{ex:gauss}, we get the following
cases of uniformly dominant convergence of $S_3(n)$:
\begin{itemize}
\item $|\alpha|<\gamma$ and $\re{(d+f-a-b-c)}>0$.
\item $|\alpha|=\gamma$, $\re{(d+f-a-b-c)}>0$ and $\re{(b+c-f)}<0$.
\item $\alpha=\gamma<0$, $\re{(b+c-f)}<-\left|\re{(d-a)}\right|$.
\end{itemize}
In particular, if the series $S_3(n)$ is balanced (that is, if $d+f=a+b+c+1$ and $\alpha=\gamma$), it converges uniformly in the following two situations:
\begin{itemize}
\item If $\gamma>0$ and $\re{(b+c-f)}<0$; 
see also \cite[Section 3]{koornwinder}. 
\item If $\gamma<0$ and $\re{(b+c-f)}<-\frac12$.
\end{itemize}
We can replace $\re{(b+c-f)}$ by $\re{(d-a-1)}$ in these two conditions.


\end{example}

\section{Proof of the main theorem}
\label{mainproof}

Here we prove Theorem \ref{th:uctheorem}. The strategy is outlined by Lemma
\ref{hpguc}.

With our notation and summation/product conventions, we may split the
hypergeometric summand $u(n,k)$ in the following ways. Firstly, we can
switch to variables in (\ref{no:hacs})--(\ref{no:tacs}) as follows:
\begin{equation} \label{eq:snkexpr}
u(n,k)={\prod_{\beta_j=0}\left({\textstyle\widetilde{a}_j+\frac{1-\alpha_j}2}\right)_{\alpha_j\,n}}
\;\xi^n\,P(n,k)\,{\prod\left({\textstyle\widehat{a}_j+\frac{1-\beta_j}2}\right)_{\beta_j\,k+\alpha_j\,n}}
\;\frac{\theta^k}{k!}.
\end{equation}
In another way, we can split the Pochhammer symbols in other way and obtain
the following expression for $u(n,k)$:
\begin{equation} \label{eq2:snkexpr}
{\prod\left({\textstyle\widetilde{a}_j+\frac{1-\alpha_j}2}\right)_{\alpha_j\,n}}\,\xi^n\,
P(n,k)\,
{\prod_{\alpha_j\neq0}\left({\textstyle\widehat{a}_j+\frac{1-\beta_j}2+\alpha_jn}\right)_{\beta_j\,k}}
\,{\prod_{\alpha_j=0}\,\left({\textstyle\widehat{a}_j+\frac{1-\beta_j}2}\right)_{\beta_j\,k}}
\;\frac{\theta^k}{k!}.
\end{equation}
Note that here the first two factors do not depend on $k$, and the last two
terms do not depend on $n$. We will use these expressions in different cases
of Lemma \ref{hpguc}.

Condition {\em (a)} of Lemma \ref{hpguc} is satisfied under the following
necessary and sufficient restrictions:
\begin{enumerate}
\item[\em (a1)] $D_0\le 0$.
 \item[\em (a2)] If $D_0=0$, then $|z_0|\le 1$.
 \item[\em (a3)] If $D_0=0$, $|z_0|=1$, then $D^*_1\le0$
and $\re{A_0}<0$.
\end{enumerate}
because for fixed general $n$ we have
\begin{eqnarray}
u(n,k)&\sim&
C(n)\,k^{\deg_kP(n,k)+\sum\left(\widehat{a}_j+\alpha_jn\right)}\;
\frac{\theta^k}{k!}\,\Gamma(k)^{\sum\beta_j}\,
{\prod\beta_j^{\beta_jk}}\\
&\sim&C(n)\,k^{D^*_1\,n+A_0-1}\,z_0^k\,\Gamma(k)^{D_0}.
\end{eqnarray}
Recall that $k!=k\,\Gamma(k)$. These conditions are general convergence
conditions for hypergeometric series; see \cite[Theorems
2.1.1--2]{specfaar}.

For condition {\em (b)} of Lemma \ref{hpguc}, we fix general $k$ and use
(\ref{eq2:snkexpr}), Corollary \ref{le:gamman}:
\begin{eqnarray} \label{as:nlimit}
u(n,k)&\sim& 
n^{\deg_nP(n,k)+\sum\left(\widetilde{a}_j+\beta_jk\right)}\,\xi^n\,
\Gamma(n)^{\sum\alpha_j}\,{\prod\alpha_j^{\alpha_j\,n}}\nonumber\\
&&\times\,H_0\,Q(k)\,\frac{\theta^k}{k!}\,{\prod_{\alpha_j\neq0}\alpha_j^{\,\beta_jk}}\,
{\prod_{\alpha_j=0}\,\left({\textstyle\widehat{a}_j+\frac{1-\beta_j}2}\right)_{\beta_j\,k}}.
\end{eqnarray}
Here $H_0$ and $Q(k)$ are the same as in (\ref{no:h0}), (\ref{no:qk}).
The first line of the right-hand side can be rewritten as 
$n^{D_1k+A^*_\infty}\,\zeta_0^n\,\Gamma(n)^{D^*_0}$. The second line is
independent of $n$. For the existence of the termwise limit we first check
whether $u(n,k)$ is bounded as $n\to\infty$, and whether the limit
$\lim_{n\to\infty} u(n,0)$ exists:
\begin{itemize}
\item[\em (b1)] $D^*_0\le 0$.
 \item[\em (b2)] If $D^*_0=0$, then $|\zeta_0|\le 1$.
 \item[\em (b3)] If $D^*_0=0$, $|\zeta_0|=1$, then $D_1\le 0$ and
$\re{A^*_\infty}\le 0$.
 \item[\em (b4)] If $D^*_0=0$, $|\zeta_0|=1$,
$\re{A^*_\infty}=0$, then $\zeta_0=1$ and $A^*_\infty=0$.
\end{itemize}
Under these conditions the termwise limit $\lim_{n\to\infty}\bS(n)$ is the
zero series if $D^*_0<0$, $|\zeta_0|<1$ or $\re{A^*_\infty}<0$. Otherwise
condition {\em (b4)} applies. Then 
the termwise limit is $H_0\,Q(0)$ if $D_1<0$, and it is equal to (\ref{glimzn})
if $D_1=0$. In these cases, asymptotics (\ref{as:nlimit}) can be rewritten,
up to a constant factor, as $k^{A^*_0-1}\,z_1^k \,\Gamma(k)^{D_0-D_1}$.
Additional conditions for the convergence of the limit series are the
following:
\begin{itemize}
\item[\em (b5)] If $D^*_0=0$, $\zeta_0=1$, $D_1=0$, $A^*_\infty=0$, then
$D_0\le0$.
 \item[\em (b6)] If $D^*_0=0$, $\zeta_0=1$,
$D_1=0$, $A^*_\infty=0$, $D_0=0$, then $|z_1|\le 1$.
 \item[\em (b7)] If $D^*_0=0$, $\zeta_0=1$, $D_1=0$, $A^*_\infty=0$,
 $D_0=0$, $|z_1|=1$, then $\re{A_0^*}<0$.
\end{itemize}

Now we check condition {\em (c)} of Lemma \ref{hpguc}. We assume that
$N=N(k)$ is an integer-valued function such that $N(k)\sim C_0\,k^p$ as
$k\to\infty$, with $p>1$ and $C_0>0$ real constants. Using formula
(\ref{eq2:snkexpr}), Corollaries \ref{gammaas} and \ref{pocchasnbig} we get
the following asymptotic expression as $k\to\infty$:
\begin{eqnarray} \label{termsnn0}
u(N,k)&\sim&C_1\,N^{\sum\wtilde{a}_j}\, {\prod\alpha_j^{\alpha_jN}}\;
\Gamma(N)^{\sum\alpha_j}\,\xi^N\,P(N,k)\,{\prod_{\alpha_j\neq
0}\alpha_j^{\beta_j\,k}}\; N^{\sum_{\alpha_j\neq0}\beta_jk}\nonumber\\
&&\times\,
\exp\left(k\,\widetilde{\Phi}_\infty\!\left(\frac{k}N\right)\right)\,
k^{\sum_{\alpha_j=0}\widehat{a}_j}\,
{\prod_{\alpha_j=0}\beta_j^{\beta_j\,k}}\;
\Gamma(k)^{\sum_{\alpha_j=0}\beta_j}\,\frac{\theta^k}{k!},
\end{eqnarray}
for some $C_1\in\RR$, and
\begin{equation} \label{phiifunc}
\widetilde{\Phi}_\infty(x)=\sum_{\alpha_j\neq
0}\beta_j\,\Theta\!\left(\frac{\beta_j\,x} {\alpha_j}\right).
\end{equation}
We rearrange as
\begin{eqnarray}  \label{termsnn}
u(N,k)&\sim&C_1\,\Gamma(N)^{D^*_0}\;\,\zeta_0^N\,\Gamma(k)^{D_0}
\;\left(\frac{N^k}{\Gamma(k)}\right)^{D_1}\,
z_1^k\;\exp\left(k\;\widetilde{\Phi}_\infty\!\left(\frac{k}{N}\right)\right) \nonumber \\
&&\hspace{-6pt}\times P(N,k) \;k^{\left(\sum\widetilde{a}_j\right)\,p\,+
\sum_{\alpha_j=0}\widehat{a}_j-1}.
\end{eqnarray}
We compute that, as $k\to\infty$,
\begin{eqnarray} \label{logbign}
\frac{\log |u(N,k)|}{k} &\equal& D^*_0\,\frac{N\log N-N}{k}+
\frac{N}{k}\,\log|\zeta_0|+D_0\left(\log k-1\right)\nonumber\\ &&
+D_1\left(\log\frac{N}{k}+1\right) +\log\left|z_1\right|+ o(1).
\end{eqnarray}
Note that $\log(N/k)\sim (p-1)\,\log k +O(1)$.

To investigate absolute convergence of $\sum_{k=0}^{\infty} u(N,k)$, we
first look at formula (\ref{logbign}) and use Lemma \ref{stconverge} with
$\varrho=1$. The series must converge absolutely for all relevant $N=N(k)$.
The most subtle case is when the expression in (\ref{logbign}) is $o(1)$.
Eventually we get the following list of conditions:
\begin{enumerate}
\item[\em (c1)] $D^*_0\le 0$.
 \item[\em (c2)] If $D^*_0=0$, then $|\zeta_0|\le 1$.
 \item[\em (c3)] If $D^*_0=0$, $|\zeta_0|=1$, then $D_0\le0$ and
$D_1\le 0$.
 \item[\em (c4)] If $D^*_0=0$, $|\zeta_0|=1$, $D_0=0$,
$D_1=0$, then $|z_1|\le 1$.
 \item[\em (c5)] If $D^*_0=0$, $|\zeta_0|=1$, $D_0=0$,
$D_1=0$, $|z_1|=1$, then $\re{A^*_\infty}\le 0$ and one of the following
conditions holds:
\begin{itemize}
\item[\em (c5A)] $\widetilde{\Phi}_\infty(x)\equiv 0$, and $\re{A_1^*}<0$.
 \item[\em (c5B)] $\widetilde{\Phi}_\infty(x)\equiv 0$, $\re{A_1^*}=0$, and either
$\deg_{\{n,k\}} P(n,k)>\deg_{n}P(n,k)+\deg_kQ(k)$ or $\re{A^*_\infty}<0$.
 \item[\em (c5C)] $\widetilde{\Phi}_\infty(x)=v_mx^m+O(x^{m+1})$ around $x=0$ for
some positive integer $m$ and negative real $v_m$, and
$\mbox{Re}\;{\psi_\infty\!\left(\frac{m+1}m\right)}<0$.
\end{itemize}
\end{enumerate}
Here we comment the case when the expression in (\ref{logbign}) is $o(1)$ as
$k\to\infty$. Formula (\ref{termsnn}) becomes then, for general $N(k)$ by
the first two parts of Lemma \ref{lem:phip},
\begin{eqnarray} \label{termsnnsi}
u(N,k)\sim \widetilde{C}_1\,
\exp\left(k\;\widetilde{\Phi}_\infty\!\left(\frac{k}{N}\right)\right)\,
k^{\,\psi_\infty\left(p\right)-1},
\end{eqnarray}
for some $\widetilde{C}_1\in\RR$. To have convergence for large $p$, we must
have $\re{A^*_\infty}\le 0$. If $\widetilde{\Phi}_\infty(x)\equiv0$ we must
have $\re{\psi_\infty(p)}<0$ for all $p\in(1,\infty)$. By part {\em (ii)} of
Lemma \ref{lem:psip}, the real part of $\psi_\infty(p)$ approaches its
supremum with $p\mapsto 1$. The condition $\re{\psi_\infty(p)}<0$ is ensured
in Case {\em (c5A)}. The Case {\em (c5B)} occurs when the supremum is not
achieved inside the interval $(1,\infty)$. If
$\widetilde{\Phi}_\infty(x)\not\equiv0$, then the exponential factor in
(\ref{termsnnsi}) is asymptotic to
\begin{equation}
\exp\left( \frac{v_m}{C_0^m}\;k^{1-(p-1)\,m} \right).
\end{equation}
For $p\ge \frac{m+1}m$ then the exponential factor is asymptotically a
constant. Then we must have $\re{\psi_\infty(p)}<0$ for all $p\in
\left[\frac{m+1}m,\infty\right)$; by part {\em (iii)} of Lemma
\ref{lem:psip} we have to check the value
$\re{\psi_\infty\!\left(\frac{m+1}m}\right)$. If $p<\frac{m+1}m$ then the
exponential factor determines convergence; the condition on $v_m$ follows
from Lemma \ref{stconverge} with $\varrho=1-(p-1)\,m$.

Now we check condition {\em (d)} of Lemma \ref{hpguc}. We assume that
$N(k)\sim C_0\,k^p$, where $p\in (0,1)$ and $C_0>0$ are real constants.
Using formula (\ref{eq:snkexpr}), Corollary \ref{gammaas} and Lemma
\ref{pocchasnl}, we get the following asymptotic expression as $k\to\infty$:
\begin{eqnarray} \label{termsnn0}
u(N,k)&\sim&C_1\,N^{\sum_{\beta_j=0}\wtilde{a}_j}\,
{\prod_{\beta_j=0}\alpha_j^{\alpha_jN}}\;
\Gamma(N)^{\sum_{\beta_j=0}\alpha_j}\,\xi^N\, {\prod_{\beta_j\neq
0}\beta_j^{\alpha_jN}}\; {\prod\beta_j^{\beta_j\,k}}\;
k^{\sum_{\beta_j\neq0}\alpha_jN}\nonumber\\
&&\times\, \exp\left(N\,\widetilde{\Phi}_0\!\left(\frac{N}k\right)\right)\,
P(N,k)\, k^{\sum\widehat{a}_j}\,
\Gamma(k)^{\sum\beta_j}\,\frac{\theta^k}{k!},
\end{eqnarray}
where $C_1\in\RR$, and
\begin{equation} \label{phi0func}
\widetilde{\Phi}_0(x)=\sum_{\beta_j\neq
0}\alpha_j\,\Theta\!\left(\frac{\alpha_j\,x} {\beta_j}\right).
\end{equation}
We rearrange as
\begin{eqnarray} \label{nkisok}
u(N,k)&\!\!\!\sim\!\!\!&C_0\;\Gamma(k)^{s-r-1}\;z_0^k\;
\Gamma(N)^{S-R}\;\left(
\frac{k^N}{\Gamma(N)}\right)^{\widetilde{S}-\widetilde{R}}\,\zeta_1^N\,
\exp\!\left(N\,\widetilde{\Phi}_0\!\left(\frac{N}{k}\right)\!\right)\nonumber\\
&& \times P(n,k)\,k^{\sum\widehat{a}_j-\sum\widehat{c}_j+
\left(\sum_{\beta_j=0}\widetilde{a}_j
-\sum_{\delta_j=0}\widetilde{c}_j\right)p-1},
\end{eqnarray}
We compute that
\begin{eqnarray} \label{lognkok}
\frac{\log |s(N,k)|}{k}&\equal&D_0\,(\log k-1)+
\log|z_0|+D^*_0\,\frac{N}{k}\big(p\,\log k\!-\!1\big)\nonumber\\
&&+D^*_1\,\frac{N}k\big((1-p)\log k+1\big)
+\frac{N}{k}\,\log\left|\zeta_1\right|+o\left(k^{-1+p}\right).
\end{eqnarray}
The last two expressions can be conveniently compared with
(\ref{termsnn})--(\ref{logbign}). Currently, $k\gg N$. Like in the previous
case, first we consider formula (\ref{lognkok}) and use Lemma
\ref{stconverge} with $\varrho=p$. 
We get a similiar set of conditions:
\begin{enumerate}
\item[\em (d1)] $D_0\le0$.
 \item[\em (d2)] If $D_0=0$ then $|z_0|\le 1$.
 \item[\em (d3)] If $D_0=0$, $|z_0|=1$, then $D^*_0\le0$ and $D^*_1\le0$.
 \item[\em (d4)] If $D_0=0$, $|z_0|=1$, $D^*_0=0$, $D^*_1=0$,
then $|\zeta_1|\le 1$.
 \item[\em (d5)] If $D_0=0$, $|z_0|=1$, $D^*_0=0$, $D^*_1=0$,
$|\zeta_1|=1$, then $\re{A_0}<0$ and one of the following conditions holds:
\begin{itemize}
\item[\em (d5A)] $\widetilde{\Phi}_0(x)\equiv 0$ and $\re{A_1}\le 0$.
 \item[\em (d5B)]
$\widetilde{\Phi}_0(x)=v_mx^m+O(x^{m+1})$ around $x=0$ for some positive
integer $m$ and negative real $v_m$, and
$\re{\psi_0\!\left(\frac{m}{m+1}\right)}<0$.
\end{itemize}
\end{enumerate}
In condition {\em (d5)}, we may consider possibilities for $\re{A_0}=0$, but
this is unnecessary because of condition {\em (a3)}. In condition {\em
(d5A)}, the case $\re{A_1}=0$ ought to be supplemented by conditions that
$\re{\psi_0(p)}\neq 0$ for all $p<1$; but this is obsolete, since if the
linear slope of $\re{\psi_0(p)}$ immediately to the left of $p=1$ is zero,
then the supremum is approached with $p\to 0$ by part {\em (iii)} of Lemma
\ref{lem:psip}.

The case when the expression in (\ref{lognkok}) is $o(k^{-1+p})$ is similar
to the consideration of $o(1)$ in (\ref{logbign}). Formula (\ref{nkisok})
becomes then, for general $N(k)$,
\begin{eqnarray} \label{termsnns}
u(N,k)\sim \widetilde{C}_1\,
\exp\left(N\;\widetilde{\Phi}_0\!\left(\frac{N}{k}\right)\right)\,
k^{\psi_0\left(p\right)-1},
\end{eqnarray}
for some $\widetilde{C}_1\in\RR$.
If $\widetilde{\Phi}_0(x)\equiv0$ we must have $\re{\psi_0(p)}<0$ for all
$p\in(0,1)$. By part {\em (ii)} of Lemma \ref{lem:psip}, we have to check
the behavior of $\psi_0(p)$ near the end-points $p=0$ and $p=1$. If
$\widetilde{\Phi}_0(x)\not\equiv0$, then the exponential factor in
(\ref{termsnns}) is asymptotic to
\begin{equation}
\exp\left( v_mC_0^{m+1}\,k^{p-(1-p)\,m} \right).
\end{equation}
For $p\le \frac{m}{m+1}$ then the exponential factor is asymptotically a
constant. Then we must have $\re{\psi_0(p)}<0$ for all $p\in
\left(0,\frac{m}{m+1}\right]$; by part {\em (iii)} of Lemma \ref{lem:psip}
we have to check the values ${\psi_0\!\left(\frac{m}{m+1}\right)}$ and
${\psi_0\!\left(0\right)}$. If $p>\frac{m}{m+1}$ then the exponential factor
determines convergence; the condition on $v_m$ follows from Lemma
\ref{stconverge} with $\varrho=p-(1-p)\,m$.

It remains to check condition {\em (e)} of Lemma \ref{hpguc}. Let us define
the family of functions:
\begin{eqnarray} \label{philfunc} \widetilde{\Phi}_\lambda(x) =
\sum_{\alpha_j\lambda+\beta_j\neq0}\alpha_j\,
\Theta\!\left(\frac{\alpha_j\,x}{\alpha_j\,\lambda+\beta_j}\right).
\end{eqnarray}
We split condition {\em (e)} into two cases:
\begin{itemize}
\item[($\star$)] $N(k)=tk+\omega(k)$ with real positive $t\not\in\Omega$,
and either $\omega(k)=O(1)$ or $\omega(k)\sim C_0k^p$ for some real
$p\in(0,1)$ and $C_0$.
 \item[($\star\star$)] $N(k)=tk+\omega(k)$ with real positive $t\in\Omega$,
and either $\omega(k)=O(1)$ or $\omega(k)\sim C_0k^p$ for some real
$p\in(0,1)$ and $C_0$.
\end{itemize}
Recall that $\Omega$ is defined in (\ref{omegadef}).

For case ($\star$) we use formula (\ref{eq4:properterm}) 
and Lemma \ref{pocchasnl} to derive the following asymptotic expression as
$k\to\infty$:
\begin{eqnarray} \label{no:termsnn0}
u(N,k)&\sim&C_1\,\Gamma(k)^{\sum\left(\beta_j+\alpha_jt\right)}
\;\frac{\theta^k}{k!}\;\xi^N\,P(N,k)\,k^{\sum\left(\widehat{a}_j+\frac12\alpha_jt\right)+
\sum_{\beta_j=0}\left(\wtilde{a}_j+\frac{t-1}2\alpha_j\right)} \nonumber\\
&&\times\,
k^{\sum\alpha_j\,\omega(k)}\,{\prod(\beta_j+\alpha_jt)^{\beta_jk+\alpha_jN}}
\,\exp\!\left(\!\omega(k)\,\widetilde{\Phi}_t\!\left(\frac{\omega(k)}k\right)\right).
\end{eqnarray}
for some $C_1\in\CC$. We arrange as follows:
\begin{eqnarray} \label{termsnn0}
|u(N,k)|&\!\sim\!& \wtilde{C}_1\,\Gamma(k)^{\sum\alpha_jt+\sum\beta_j}
\;\frac{|P(N,k)|}{k!}\,\exp\left(\omega(k)\,\widetilde{\Phi}_t\!\left(\frac{\omega(k)}k\right)\right)\nonumber\\
&&\hspace{-10pt}\times\left(|\theta||\xi|^t\,{\prod|\beta_j+\alpha_jt|^{\beta_j+\alpha_jt}}\right)^k
\left(|\xi|\,k^{\sum\alpha_j}\,{\prod|\beta_j+\alpha_jt|^{\alpha_j}}\right)^{\omega(k)}\; \nonumber\\
&&\hspace{-10pt} \times\,k^{\scre{\sum\widehat{a}_j+
\sum_{\beta_j=0}\wtilde{a}_j}+\frac{t-1}2\sum\alpha_j+
\frac{1}2\sum_{\beta_j\neq0}\alpha_j}.
\end{eqnarray}
for some $\widetilde{C}_1\in\RR$. Using (\ref{no:termstr}), (\ref{no:gt}),
(\ref{gtderiv}), we rewrite:
\begin{eqnarray} \label{asymptnkt}
|u(N,k)| & \sim & 
\Gamma(k)^{D^*_0t+D_0}\, 
\left(\frac{k^{\omega(k)+\frac{t-1}2}}{\exp(1)}\right)^{\!D^*_0}
\,g(t)^k\;\exp\left(\omega(k)\,\frac{g'(t)}{g(t)}\right) \nonumber\\
&&  \times
\exp\left(\omega(k)\;\widetilde{\Phi}_t\!\left(\frac{\omega(k)}{k}\right)\right)\;
k^{\mbox{\scriptsize Re}\,\psi^{\star}_t(p)\,\!-\,1}.
\end{eqnarray}
Here we set $p=0$ if $\omega(k)=O(1)$.
Recall that $\psi^{\star}_t(p)$ is a monotone nondecreasing function by
part {\em (ii)} of Lemma \ref{lem:phipst}.

We already have $D_0\le 0$ and $D^*_0\le0$ by conditions {\em (a1)} and {\em
(b1)}. Case ($\star$) gives additional conditions if $D_0=0$ and $D^*_0=0$.
Firstly, we must have $g(t)\le 1$ for all positive
$t\in\RR\setminus\Omega$. 
If this is the case, and $g(t_0)=1$ for some positive
$t_0\in\RR\setminus\Omega$, then $g'(t_0)=0$. Indeed, $g(t_0)\neq 0$ would
imply $g(t_1)>1$ for some $t_1\in\RR\setminus\Omega$ in a neighborhood of
$t_0$. Therefore we may ignore the exponential factor 
with $g'(t)$. At these points $t_0$ we have to consider the last two terms
in (\ref{asymptnkt}). Eventually we get the following conditions for the
case ($\star$):
\begin{itemize}
\item[\em (e1)] If $D_0=0$, $D^*_0=0$, then $g(t)\le 1$ for all positive
$t\in\RR\setminus\Omega$.
 \item[\em (e2)] If $D_0=0$, $D^*_0=0$, and $g(t)=1$ for
some positive $t\in\RR\setminus\Omega$, then for any
$t_0\in\RR\setminus\Omega$ where $g(t_0)=1$, we must have
$\re{\psi^{\star}_{t_0}(0)}<0$ and one of the following two conditions
satisfied:
\begin{itemize}
\item[\em (e2A)] $\widetilde{\Phi}_{t_0}(x)\equiv 0$, and
$\re{\psi^{\star}_{t_0}(1)}\le 0$.
 \item[\em (e2B)] $\widetilde{\Phi}_{t_0}(x)=v_mx^m+O(x^{m+1})$ around $x=0$, where $m$ is a
positive {\bf odd} integer, $v_m<0$, and
$\re{\psi^{\star}_{t_0}\!\left(\frac{m}{m+1}\right)}<0$.
\end{itemize}
\end{itemize}
Here we comment the situations when condition {\em (e2)} applies. We have
$\re{\psi^{\star}_{t_0}(0)}<0$ because the power of $k$ in (\ref{asymptnkt})
determines the convergence when $\omega(k)=O(1)$. If
$\widetilde{\Phi}_t(x)\equiv0$, we must have $\re{\psi^{\star}_t(p)}<0$ for
all $p\in[0,1)$. By part {\em (iii)} of Lemma \ref{lem:psipstr}, it is
enough to have \mbox{$\re{\psi^{\star}_{t_0}(1)}\le 0$}. If
$\widetilde{\Phi}_t(x)\not\equiv0$, then the exponential factor is
asymptotic to $\exp\left(v_m\,C_0^{m+1}\,k^{p-(1-p)\,m}\right)$; it is
relevant when $p\in\left(\frac{m}{m+1},1\right)$. If $m$ is even, the
exponential factor is unbounded either when $C_0>0$ or when $C_0<0$. Hence
$m$ must be odd. Then Lemma \ref{stconverge} with $\varrho=p-(1-p)\,m$ gives
the restriction $v_m<0$. The power of $k$ factor must be restricted for
$p\in\left[\frac{m}{m+1},1\right]$. By part {\em (iii)} of Lemma
\ref{lem:psipstr}, it is enough to have
\mbox{$\re{\psi^{\star}_{t_0}\!\left(\frac{m}{m+1}\right)}\le 0$}.

Now we consider the case ($\star\star$), with $t\in\Omega$. Formula
(\ref{no:termsnn0}) should be modified as follows:
\begin{itemize}
\item The sums and products should be supplemented by the condition
\mbox{$\beta_j+\alpha_j\,t\neq 0$}. This is unnecessary for the sums in the
power of $\Gamma(k)$, and eventually in some products (since $0^0=1$). Note
that these conditions are already indicated in definition (\ref{philfunc})
of $\widetilde{\Phi}_t(x)$.
 \item By Lemma \ref{pochhntk}, we have to append
\begin{equation}
|\omega(k)|^{\sum_{\alpha_jt+\beta_j=0}\left(b_j-\frac12\right)}
\left(\!\sqrt{|\omega(k)|}\;\Gamma(|\omega(k)|)\right)^{\varepsilon\sum_{\alpha_jt+\beta_j=0}\alpha_j}
{\prod_{\alpha_jt+\beta_j=0}(\varepsilon\alpha_j)^{\alpha_j\omega(k)}}.
\end{equation}
\end{itemize}
With these modifications, asymptotic expression (\ref{asymptnkt}) can be
written eventually as
\begin{eqnarray} \label{nn:asymptnkt}
|u(N,k)| & \sim & 
\Gamma(k)^{D^*_0\,t+D_0}\! 
\left(\!\sqrt{|\omega(k)|}\,
\Gamma(|\omega(k)|)\right)^{\varepsilon\sum_{\alpha_jt+\beta_j=0}\alpha_j}
\;g(t)^k\nonumber\\
&&\times \left(|\xi|
{\prod_{\alpha_jt+\beta_j\neq0}\!|\beta_j+\alpha_jt|^{\alpha_j}}\!
{\prod_{\alpha_jt+\beta_j=0}\!|\alpha_j|^{\alpha_j}}\!\right)^{\!\omega(k)}\!
\exp\!\left(\omega(k)\,\widetilde{\Phi}_t\!\left(\frac{\omega(k)}{k}\right)\!\right)\nonumber\\
&&\times\;k^{D^*_0\left(\omega(k)+\frac{t-1}2\right)
-\sum_{\alpha_jt+\beta_j=0}\left(\alpha_j\omega(k)+\frac12p\beta_j\right)}
\;k^{\mbox{\scriptsize Re}\,\psi^{\star}_t(p)\,-\,1}.
\end{eqnarray}
Here we set $p=0$ if $\omega(k)=O(1)$.

As in the case ($\star$), there are extra conditions only if $D_0=0$ and
$D^*_0=0$. Then we have:
\begin{eqnarray} \label{eq:unkstar}
\log |u(N,k)| &=&
\sum_{\alpha_jt+\beta_j=0}\alpha_j\;\omega(k)\,\big( \log|\omega(k)|-\log k-1\big) \vspace{3pt}\nonumber\\
&&+\,k\log g(t)+O(w(k)+\log k).
\end{eqnarray}
In general, the dominant term is $k\log\,g(t)$; hence we must have $g(t)\le
1$.

Suppose that $g(t_0)=1$ for some $t_0\in\Omega$. If
$\Sigma_{\alpha_jt_0+\beta_j=0}\,\alpha_j\neq0$, then the first term in (\ref{eq:unkstar})
approaches $+\infty$ for those $\omega(k)\sim C_0k^p$ with $p$ close to $1$
and with $C_0>0$ or $C_0<0$ depending on the sign of $\Sigma_{\alpha_jt_0+\beta_j=0}\,\alpha_j$.
Hence $\Sigma_{\alpha_jt_0+\beta_j=0}\,\alpha_j=0$ for those $t_0\in\Omega$ for which
$g(t_0)=1$. Then $\Sigma_{\alpha_jt_0+\beta_j=0}\,\beta_j=0$ as well.

If $\Sigma_{\alpha_jt_0+\beta_j=0}\,\alpha_j=0$, then $g(t)$ is actually
differentiable at $t_0$ by part {\em (iii)} of Lemma \ref{gtgenlemma}. The
value of the derivative can be derived from (\ref{no:gtderiv}) or
(\ref{as:dgt}). If $D_0=0$, $D^*_0=0$, $g(t_0)=1$ and
$\Sigma_{\alpha_jt_0+\beta_j=0}\,\alpha_j=0$, we can rewrite
(\ref{nn:asymptnkt}) as follows:
\begin{eqnarray} \label{no:unkom}
|u(N,k)| \sim 
\exp\left(\omega(k)\,\frac{g'(t_0)}{g(t_0)}\right)
\exp\!\left(\omega(k)\,\widetilde{\Phi}_{t_0}\!\left(\frac{\omega(k)}{k}\right)\right)
k^{\mbox{\scriptsize Re}\,\psi^{\star}_t(p)\,-\,1}.
\end{eqnarray}
If $g'(t_0)\neq 0$, then condition {\em (e1)} is contradicted for some point
$t\in\RR\setminus\Omega$ in a neighborhood of $t_0$. Hence we may assume
$g'(t_0)=0$. Eventually we get the following conditions:
\begin{itemize}
\item[\em (e3)] If $D_0=0$, $D^*_0=0$, then $g(t)\le 1$ for all positive
$t\in\Omega$.
 \item[\em (e4)] If $D_0=0$, $D^*_0=0$, $g(t_0)=1$ for some positive
$t_0\in\Omega$, then for any $t_0\in\RR\setminus\Omega$ where $g(t_0)=1$ we
must have $\sum_{\alpha_jt_0+\beta_j=0}\,\alpha_j=0$,
$\re{\psi^{\star}_{t_0}(0)}<0$, and one of the following two conditions
satisfied:
\begin{itemize}
\item[\em (e4A)] $\widetilde{\Phi}_{t_0}(x)\equiv 0$, and
$\re{\psi^{\star}_{t_0}(1)}\le 0$.
 \item[\em (e4B)] $\widetilde{\Phi}_{t_0}(x)=v_mx^m+O(x^{m+1})$ around $x=0$, where $m$ is a
positive odd integer, $v_m<0$, and
$\re{\psi^{\star}_{t_0}\!\left(\frac{m}{m+1}\right)}<0$.
\end{itemize}
\end{itemize}
The subcases of ({\em e4}) are derived similarly as the subcases of ({\em
e2}). Compared with conditions {\em (e1)--(e2)}, we additionally have the
condition $\sum_{\alpha_jt_0+\beta_j=0}\,\alpha_j=0$ in {\em (e4)}. But this
condition is trivially satisfied in case ($\star$), so formally we may
require it in both cases. An implicit difference between cases ($\star$) and
($\star\star$) is that the functions $\psi^{\star}_{t}(p)$ and
$\widetilde{\Phi}_{t}(x)$ can be defined simpler in case ($\star$).

Before summarizing up the derived conditions, we remark that the nonzero
Taylor coefficients (\ref{taylorthet}) of $\Theta(x)$ have the same signs as
the Taylor coefficients of the rational function
$x/(1+x)=\sum_{j=1}^{\infty} (-1)^{j+1}x^j$. The corresponding coefficients
differ the positive factor $j\,(j+\!1)$. If we replace each occurrence of
$\Theta(x)$ by $x/(1+x)$ in definitions (\ref{phiifunc}), (\ref{phi0func}),
(\ref{philfunc}) of $\widetilde{\Phi}_{\infty}(x)$, $\widetilde{\Phi}_0(x)$,
$\widetilde{\Phi}_{\lambda}(x)$, respectively, we get the rational functions
$\Phi_{\infty}(x)$, $\Phi_0(x)$, $\Phi_{\lambda}(x)$ defined in
(\ref{no:phii}), (\ref{no:phi0}), (\ref{no:phil}), respectively. The Taylor
coefficients around $x=0$ of the rational functions differ by the positive
factor $j\,(j+\!1)$ from the respective coefficients of the corresponding
$\widetilde{\Phi}$-functions. Therefore we may replace in conditions {\em
(c5)}, {\em (d5)}, {\em (e2)}, {\em (e4)} the functions
$\widetilde{\Phi}_{\infty}(x)$, $\widetilde{\Phi}_0(x)$,
$\widetilde{\Phi}_{\lambda}(x)$ by the rational functions
$\Phi_{\infty}(x)$, $\Phi_0(x)$, $\Phi_{\lambda}(x)$, respectively.

Now we summarize the conditions {\em (a1)--(a3)}, {\em (b1)--(b7)}, {\em
(c1)--(c5)}, {\em (d1)--(d5)}, {\em (e1)--(e4)}. Note that
\begin{eqnarray*}
\refi{a1}\Rightarrow\refi{b5}\ \&\ \refi{d1},\qquad
\refi{a2}\Rightarrow\refi{d2},\qquad \refi{b1}\Rightarrow\refi{c1},\qquad
\refi{b2}\Rightarrow\refi{c2},\\
\refi{c4}\Rightarrow\refi{b6},\qquad \refi{a1}\ \&\ \refi{b3}
\Rightarrow\refi{c3}, \qquad \refi{a3}\ \&\ \refi{b1} \Rightarrow\refi{d3}.
\end{eqnarray*}
Therefore we may discard the conditions \refi{b5}--\refi{b6},
\refi{c1}--\refi{c3}, \refi{d1}--\refi{d3}. Because of \refi{a3}, we can
drop the restriction $\re{A_0}<0$ in \refi{d5}. Because of \refi{b3}, we can
drop the restriction $\re{A^*_\infty}\le0$ in \refi{c5}. Besides, in cases
\refi{c5A} and \refi{c5B} we can drop condition \refi{b7}, because
$\ree{A_\infty^*+A_0^*}\le\re{A_1^*}$.

We have the following correspondence between the conditions:
\begin{eqnarray*}
\refi{a1}\ \&\ \refi{b1} \Rightarrow\refi{i},\qquad
\mbox{\refi{a2}\,--\,\refi{a3}}\Leftrightarrow\refi{ii},\qquad
\mbox{\refi{b2}\,--\,\refi{b4}}\Leftrightarrow\refi{iii},\\
\refi{e1}\ \&\ \refi{e3} \Leftrightarrow\refi{iv},\qquad
 \refi{e2}\ \&\ \refi{e4} \Leftrightarrow\refi{v},\qquad 
\mbox{\refi{d4}\,--\,\refi{d5}}\Leftrightarrow\refi{vi},\\
\mbox{\refi{c4}\,--\,\refi{c5}}\ \&\ \refi{b7}\Leftrightarrow\refi{vii}.
\end{eqnarray*}
The limit $\lim_{n\to\infty}\bS(n)$ is discussed right after the conditions
\refi{b1}--\refi{b4} here above.\QED

\section{Properties of $g(t)$}
\label{gtsection}

As mentioned in remark \refi{VI} after Theorem (\ref{th:uctheorem}),
all conditions of Theorem \ref{th:uctheorem} can be determined
algorithmically. The only less straightforward part is dealing with the
function $g(t)$ in parts \refi{iv}--\refi{v}. This is significant when
$D_0=0$ and $D^*_0=0$. Some key properties of $g(t)$ are presented
in Lemma \ref{gtgenlemma}. Here we focus on finding local extremuma
of $g(t)$. At the end, a simplified version (\ref{simplegt}) of this function
is considered thoroughly.
\begin{lemma} \label{gtconds}
In the context of Sections $\ref{sec:notation}$ and $\ref{sec:enotation}$,
suppose that $D_0=0$, $D^*_0=0$, and that conditions (ii)--(iii), (vii) of
Theorem $\ref{th:uctheorem}$ hold. Then $g(t)\le 1$ for all $t>0$ if and
only if the following conditions hold:
\begin{itemize}
\item For all $t\not\in\Omega$ such that
\begin{equation} \label{ine:cond}
|\zeta_0|\,
\prod_{\alpha_j\neq0}\left|\,t+\frac{\beta_j}{\alpha_j}\,\right|^{\alpha_j}=
1 
\end{equation}
we have
\begin{equation} \label{ine:value}
\hspace{4pt} |z_1|\; 
\prod_{\alpha_j\neq0}\left|\,t+\frac{\beta_j}{\alpha_j}\,\right|^{\beta_j}
\le 1.
\end{equation}
\item For all $t\in\Omega$ such that equality $(\ref{rpartv})$ holds, we
have $g(t)\le 1$.
\end{itemize}
If these conditions are satisfied, then $g(t)=1$ are those points
$t\not\in\Omega$ where equalities in $(\ref{ine:cond})$ and
$(\ref{ine:value})$ hold, and possibly some points $t\in\Omega$ where
equality $(\ref{rpartv})$ holds.
\end{lemma}
\proof  By parts \refi{i}--\refi{ii} of Lemma \ref{gtgenlemma}, the function
$g(t)$ is continuous on $\RR$, and it is continuously differentiable on
$\RR\setminus\Omega$. We need to investigate the behavior of $g(t)$ as $t$
approaches $+\infty$, $0$ or singularities of $g'(t)$, and find local
extremuma of $g(t)$.

As $t\to0$, then $g(t)\to |z_0|$ by part \refi{iv} of Lemma
\ref{gtgenlemma}. But $|z_0|\le1$ by part \refi{ii} of Theorem
\ref{th:uctheorem}. As $t\to\infty$, then $g(t)\sim|z_1|\exp(D_1)\,
|\zeta_0|^t\,t^{D_1}$ by part \refi{iv} of Lemma \ref{gtgenlemma}. The chain
of possible restrictions $|\zeta_0|\le 1$, $D_1\le0$, $|z_1|\le 1$ is
implied by parts \refi{iii}, \refi{vii} of Theorem \ref{th:uctheorem}. By
part \refi{iv} of Lemma \ref{gtgenlemma}, genuine points of discontinuity of
$g'(t)$ are not local extremuma.

It remains to check the local extremuma at those $t>0$ where $g'(t)$ is
actually continuous. For these points, either $t\not\in\Omega$, or
$t\in\Omega$ and equality $(\ref{rpartv})$ holds. Condition (\ref{ine:cond})
is just reformulation of $g'(t)=0$, following expression (\ref{no:gtderiv}).
Recall that we assume $S=R$. Inequality (\ref{ine:value}) is equivalent to
$g(t)\le 1$ if condition (\ref{ine:cond}) is satisfied.

If $g(t)\le 1$ for all $t>0$, then the points with $g(t)=1$ are local
extremuma. If $t\in\Omega$ and $g'(t)=0$, then the quotient of the left and
right hand sides of (\ref{ine:value}) is equal to $g(t)$. \QED

\noindent Here we continue the list of observations \refi{I}--\refi{IX}
in Section  \ref{sec:theresult} with a few more remarks.
\begin{enumerate}
\item[\refi{X}] If all $\alpha_j$'s and $\gamma_j$'s are even, then condition
(\ref{ine:cond}) is actually a polynomial equation for $t$. If there are
some odd $\alpha_j$'s or $\gamma_j$'s, we can square both sides of
(\ref{ine:cond}) and get a polynomial equation for $t$ as well. We have to
find real positive roots of these equations. The numeric or algebraic roots
of the polynomial equations can be found algorithmically. On the other hand,
the equations might have inappropriately high degree. It might be useful to
have some estimates of the number and location of relevant solutions.
\item[\refi{XI}] 
The two conditions for $g(t)\le1$ 
can be formulated in a single statement, if we add the condition
$\alpha_jt+\beta_j\neq0$ to the products in (\ref{ine:cond}) and
(\ref{ine:value}), or make the convention that the both-side factors
$|t-\lambda|$ with $\lambda=t$ in these formulas cancel out if $t\in\Omega$
and equality (\ref{rpartv}) holds. The unified statement is: For all $t>0$
such that equalities (\ref{rpartv}) and (\ref{ine:cond}) hold, we must have
(\ref{ine:value}). Identification of the points $g(t)=1$ in Lemma \ref{gtconds}
can be similarly unified. From algorithmic point of view, the single equation
(\ref{ine:cond}) with simplified or cancelled-out powers of $|t-\lambda|$
determines all local extremuma.
 \item[\refi{XII}]  Let us denote $h(t)=g'(t)/g(t)$. Using formula (\ref{gtderiv}) we
derive
\begin{equation}
h'(t)=
\sum \frac{\alpha_j^2}{\alpha_jt+\beta_j}.
\end{equation}
If we compute the zeroes and poles of this rational function, and (signs of)
values of $h(t)$ there, we can determine intervals where zeroes of $h(t)$
lie. Since $g'(t)$ has the same sign as $h(t)$ for any $t\not\in\Omega$,
those are also intervals for the zeroes of $g'(t)$, or extremuma of $g(t)$.
 \item[\refi{XIII}]  Lemma \ref{gtconds} formally holds in the case when $g(t)$ is a
constant function as well. Of course, in that case condition {\em (iv)} of
Theorem \ref{th:uctheorem} is straightforward to handle.
 \item[\refi{XIV}]  The second paragraph in the proof of Lemma \ref{gtconds} shows that
part {\em (iv)} of Theorem \ref{th:uctheorem} implies $|z_1|\le 1$ if
$|\zeta_0|=1$, $D_1=0$ (and $D_0^*=0$). Consequently, one may simplify 
part {\em (vii)} of Theorem \ref{th:uctheorem} by starting 
``{\em If $|\zeta_0|=1$, $D_1=0$ and $|z_1|=1$, then ...}", 
and dropping all further conditions on $z_1$. Similarly, because of the asymptotics
in (\ref{eq:ldgt0}),  part {\em (iv)} of Theorem \ref{th:uctheorem} implies $|\zeta_1|\le 1$ if 
$|z_0|=1$, $D^*_1=0$ (and $D_0^*=0$).  Hence, part {\em (vi)} of Theorem \ref{th:uctheorem} 
can be simplified by starting ``{\em If $|z_0|=1$, $D^*_1=0$ and $|\zeta_1|=1$, then ...}",
and dropping all further conditions on $\zeta_1$. But from computational point of view, it is
convenient to use formulation of Theorem \ref{th:uctheorem} so to handle the behavior
of $g(t)$ as $t\to 0$ and $t\to\infty$ automatically.
\end{enumerate}
\renewcommand{\labelenumi}{(\roman{enumi})}
In the rest of this Section, we explicitly consider a simple case of the
$g(t)$-function:
\begin{equation} \label{simplegt}
\widehat{g}(t)=\frac{|\alpha\,t+1|^{\alpha t+1}\;|\gamma\,t|^{\gamma t}}
{|\alpha\,t|^{\alpha\,t}\;|\gamma\,t+1|^{\gamma t+1}},
\end{equation}
This case naturally occurs with sequences of hypergeometric functions of the
form
\begin{equation}
\hpg{s}{r}{a_1+\alpha\,n,\,a_2,\,\ldots,\,a_s}{c_1+\gamma\,n,\,c_2,\,\ldots,\,c_r}{\;z}.
\end{equation}
There may be more upper and lower parameters dependant on $n$, if they
cancel each other out in the expression of $g(t)$. We saw the same function in
Example \ref{ex:gauss}; see (\ref{eq:ex71gt}). Knowledge of the function
$\widehat{g}(t)$ may help to arrive at effective estimates for more
complicated functions $g(t)$, by splitting them into a product of
$\widehat{g}(t)$'s.


In the following Lemma, we present basic properties of $\widehat{g}(t)$. We
assume here that $\gamma>0$, but allow $t$ to be both positive and negative.
If $\gamma<0$ in (\ref{simplegt}), then Lemma \ref{gtsplemma} can be applied
by considering $\gamma\mapsto-\gamma$, $t\mapsto-t$, $\alpha\mapsto-\alpha$,
so that $\gamma>0$ and $t<0$.
\begin{lemma} \label{gtsplemma}
Assume that $\alpha$, $\gamma$ are integers, $\alpha\neq\gamma$ and
$\gamma>0$. 
\begin{enumerate}
\item The function $\widehat{g}(t)$ is continuous on the whole real axis,
and is differentiable everywhere  except the points
$x\in\{0,-1/\alpha,-1/\gamma\}$. These three points are not local extremuma.
 \item $\widehat{g}(0)=1$, and $\lim_{t\to\pm\infty} \widehat{g}(t)=|\alpha/\gamma|$.
 \item $\sup_{t>0} \widehat{g}(t)=\max(1,|\alpha/\gamma|)$.
 \item The global supremum of $\widehat{g}(t)$ is
achieved for a negative $t$, and it is the only local extrema which
satisfies $\widehat{g}(t)>1$ and $\widehat{g}(t)>|\alpha/\gamma|$.
\end{enumerate}
\end{lemma}
\proof The first part follows from parts \refi{i}--\refi{iii} of Lemma
\ref{gtgenlemma}. The value $\widehat{g}(0)$ is trivial. We have
\[
\lim_{t\to\infty} \widehat{g}(t)=\lim_{t\to\infty}
\frac{|\alpha\,t+1|}{|\gamma\,t+1|} \,\left| 1+\frac{1}{\alpha\,t}
\right|^{\alpha\,t} \Big/ \left| 1+\frac{1}{\gamma\,t} \right|^{\gamma\,t} =
\left| \frac{\alpha}{\gamma} \right| ,
\]
and similarly for  $\lim_{t\to-\infty} \widehat{g}(t)$.

Let us consider
\begin{equation} \label{ghtderiv}
h(t):=\frac{\widehat{g}{\,}'(t)}{\widehat{g}(t)} =
\alpha\,\log|1+\alpha\,t|-\alpha\,\log|\alpha\,t|
+\gamma\,\log|\gamma\,t|-\gamma\,\log|1+\gamma\,t|.
\end{equation}
The local extremuma of $\widehat{g}(t)$ are determined by $h(t)=0$. We have:
\begin{equation}
h'(t) = \frac{\gamma-\alpha}{t\,(1+\alpha\,t)\,(1+\gamma\,t)},
\end{equation}
We conclude that $h(t)$ and $\widehat{g}'(t)$ are monotone on the intervals
separated by points $0$, $-1/\alpha$ and $-1/\gamma$. Here are some relevant
limits:
\begin{eqnarray*}
\lim_{t\to-1/\gamma} h(t) = \infty, & & \lim_{t\to\pm\infty} h(t) = 0,\\
\lim_{t\to-1/\alpha} h(t) = \left\{ \begin{array}{rl} -\infty, & \mbox{if }
\alpha>0,\\ \infty, & \mbox{if } \alpha<0.
\end{array} \right. & & \lim_{t\to 0} h(t) =
\left\{ \begin{array}{rl} \infty, & \mbox{if } \alpha>\gamma,\\
-\infty, & \mbox{if } \gamma>\alpha. \end{array} \right.
\end{eqnarray*}
We distinguish the following cases:
\begin{figure} 
\centering \epsfig{file=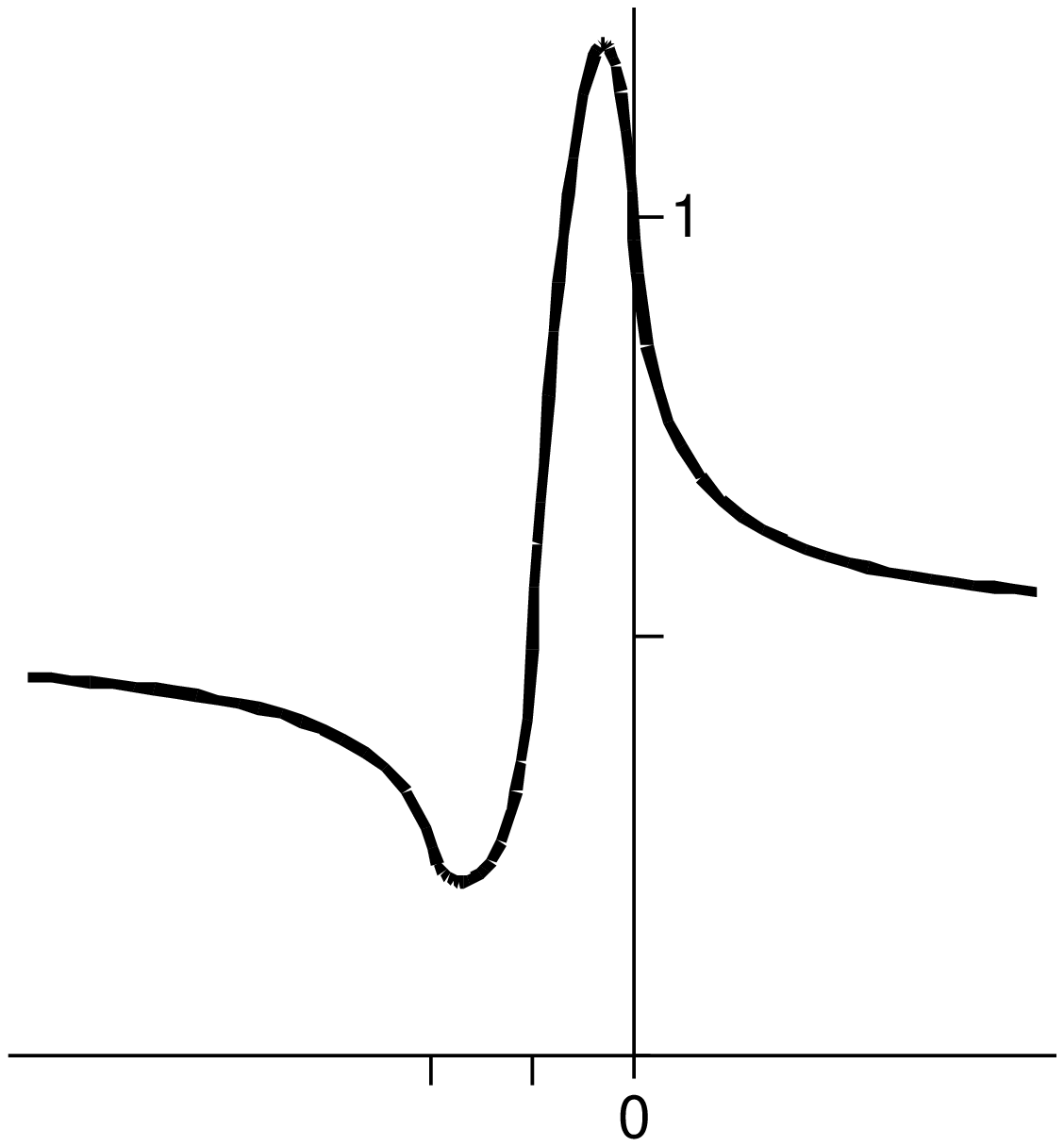, width=1.6in}  \hspace{-5pt}
\epsfig{file=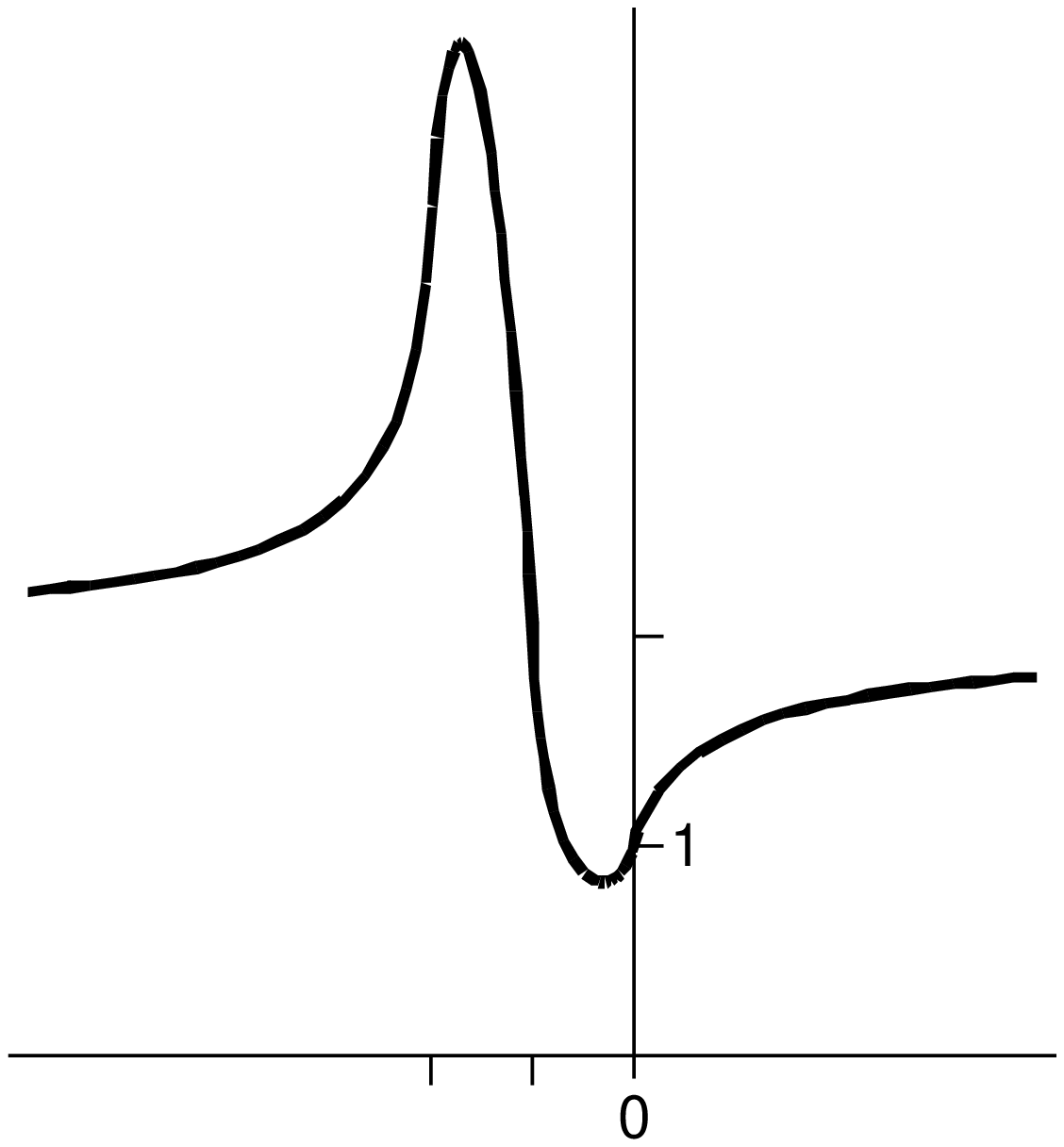, width=1.6in} \hspace{-5pt}
\epsfig{file=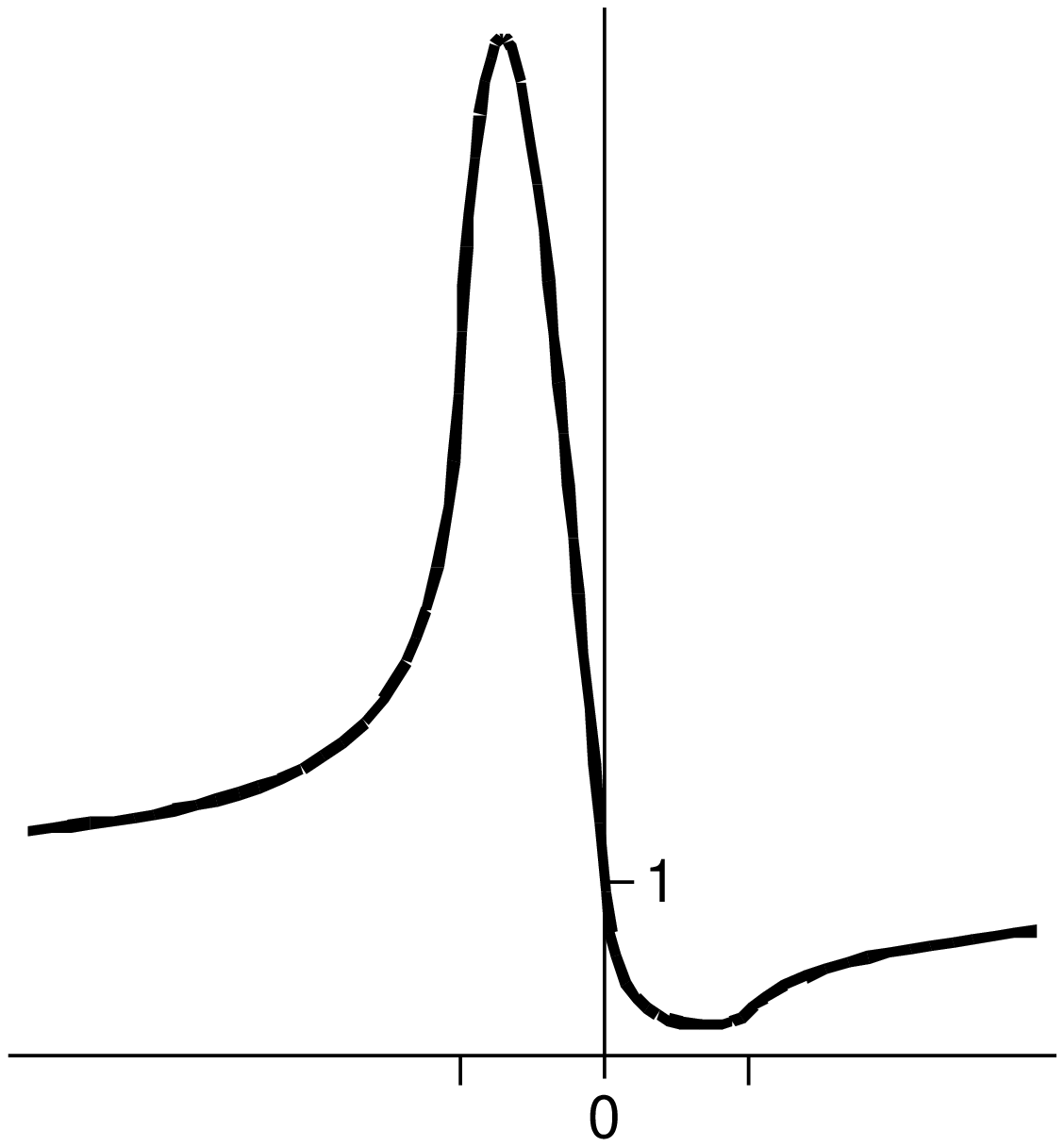, width=1.6in} \hspace{-5pt}
\epsfig{file=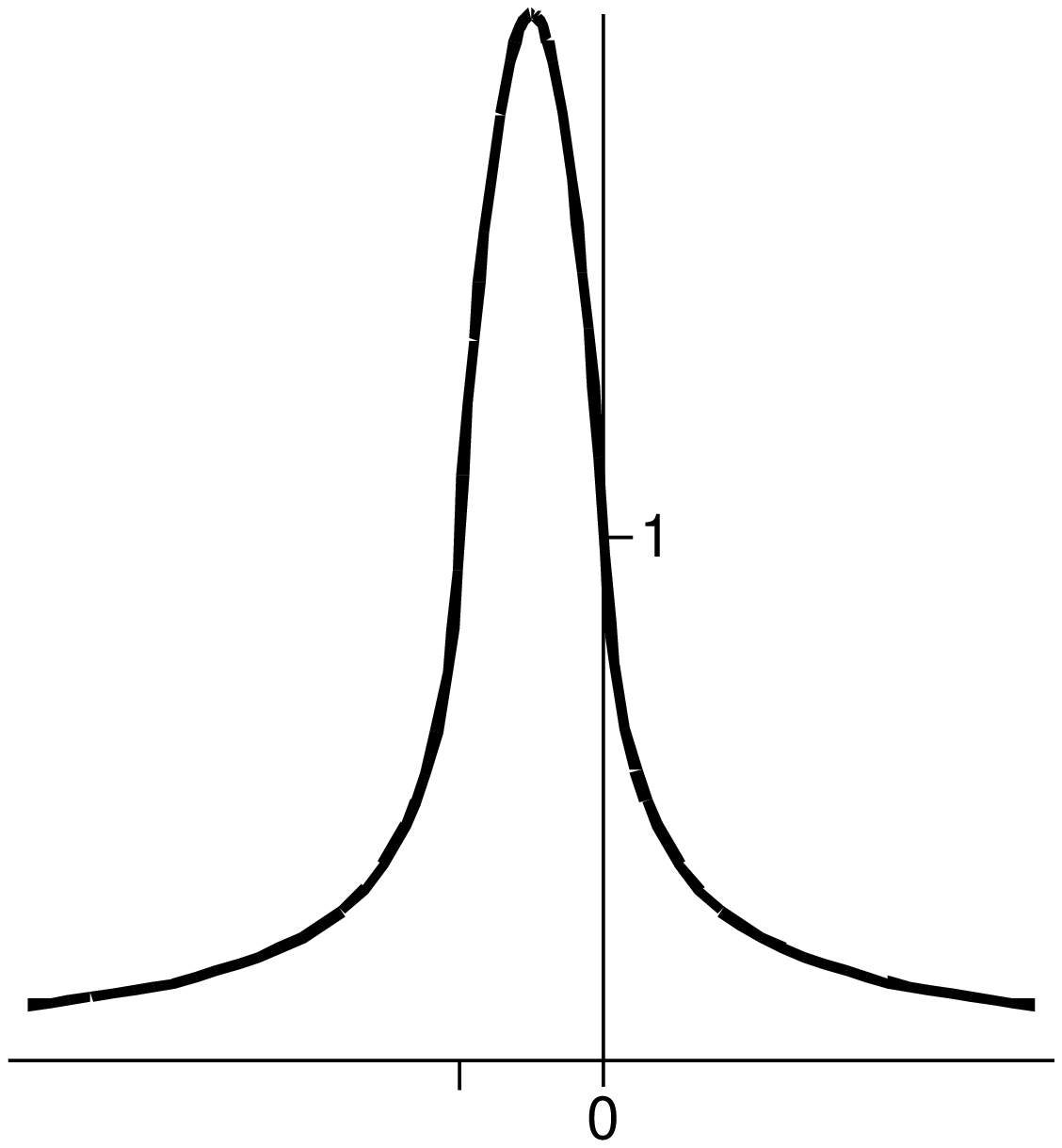, width=1.6in}
\caption{The function $\widehat{g}(t)$ for various $\alpha/\gamma$}
 \label{fig:cycle78}
\end{figure}
\begin{itemize}
\item If $0<\alpha<\gamma$, then $\widehat{g}(t)$ has a local maximum on the
interval $(-1/\gamma,0)$, which is greater than
$\widehat{g}(0)=1>\alpha/\gamma$. There is a local minimum on
$(-1/\alpha,-1/\gamma)$, which is less than
$\widehat{g}(-\infty)=\alpha/\gamma<1$. For positive $t$, the function
$\widehat{g}(t)$ decreases from 1 to $\alpha/\gamma$. See the first graph in
Figure \ref{fig:cycle78}.
 \item If $0<\gamma<\alpha$, then $\widehat{g}(t)$ has a local maximum
on the interval $(-1/\gamma,-1/\alpha)$, which is greater than
$\widehat{g}(-\infty)=\alpha/\gamma>1$. There is a local minimum on
$(-1/\alpha,0)$, which is less than $\widehat{g}(0)=1<\alpha/\gamma$. For
positive $t$, function $\widehat{g}(t)$ increases from 1 to $\alpha/\gamma$.
See the second graph in Figure \ref{fig:cycle78}.
 \item If $\alpha<0$, then $\widehat{g}(t)$ has a local maximum on the interval
$(-1/\gamma,0)$, which is greater than
$\widehat{g}(-\infty)=|\alpha/\gamma|$ and $\widehat{g}(0)=1$. There is a
local minimum on $(0,-1/\alpha)$, which is less than $1$ and
$|\alpha/\gamma|$. The supremum of $\widehat{g}(t)$ over positive $t$ is
achieved as $t\to 0$ or $t\to\infty$. See the third graph in Figure
\ref{fig:cycle78}.
 \item If $\alpha=0$, then $g(t)$ has a
local maximum on the interval $(-1/\gamma,0)$, which is greater than
$\widehat{g}(0)=1$. There are no other extremuma in this case. For positive
$t$, the function $\widehat{g}(t)$ decreases from 1 to 0. See the last graph
in Figure \ref{fig:cycle78}.
\end{itemize}
This analysis proves parts \refi{iii}--\refi{iv} of the Lemma. \QED

\begin{conclude}
Suppose that $\alpha\neq\gamma$. If $\gamma>0$, then the supremum of
$\widehat{g}(t)$ over $t>0$ is achieved either as $t\to 0$ or $t\to\infty$.
If $\gamma<0$, then the supremum of
$\widehat{g}(t)$ over $t>0$ is achieved for some 
$t\in\left(0,\frac1{|\gamma|}\right)$.
\end{conclude}
\proof If $\gamma>0$, we use parts \refi{ii}--\refi{iii} of Lemma
\ref{gtsplemma}. If $\gamma<0$ then we apply Lemma \ref{gtsplemma} after
changing the signs $\gamma\mapsto-\gamma$, $t\mapsto-t$,
$\alpha\mapsto-\alpha$.\QED

\noindent To estimate how high is the maximum of $\widehat{g}(t)$ over those
$t$ with $\gamma\,t<0$, we need this Lemma.

\begin{lemma} \label{ytlemma}
Suppose that $x\ge 1$. The equation
\begin{eqnarray} \label{eq4:ybigx}
\frac{y^x}{x^x} = \frac{(y+1)^{x-1}}{(x-1)^{x-1}}
\end{eqnarray}
has a unique root $y$ such that $y\ge 1$.

Let $y(x)$ denote the unique root as a function of $x$. Asymptotically,
\begin{equation} \label{gasympt}
y(x) \sim \tau\,x-\frac{\tau+1}{2}
-\frac{(\tau+1)\,(\tau-2)}{24\,\tau}\,\frac{1}{x} +\ldots, \qquad
\mbox{as } x\to\infty,
\end{equation}
where $\tau$ is the real solution of $\log(\tau)=1+1/\tau$:
\begin{equation} \label{tauval}
\tau\approx 3.59112147666862213664922292574163484210\ldots.
\end{equation}
For $x\ge 1$ we have
\begin{equation} \label{ineq4big}
\tau\,(x-1)+1<y(x)<\tau\,(x-1)+\frac{\tau-1}{2}.
\end{equation}
\end{lemma}
\proof Let us consider the logarithm of the ratio of both sides of
(\ref{eq4:ybigx}):
\begin{equation} \label{eq:psixy}
\Psi(x,y)=x\,\log y -x\,\log x-(x\!-\!1)\,\log(y\!+\!1)+
(x\!-\!1)\,\log(x\!-\!1).
\end{equation}
For fixed $x\ge 1$, we have to find solutions of $\Psi(x,y)=0$ with $y\ge
1$. We have: 
\begin{equation}
\frac{\partial\Psi(x,y)}{\partial y}=\frac{y+x}{y(y+1)}.
\end{equation}
Hence, as a function of $y$, $\Psi(x,y)$ is continuous increasing function
on the interval $(1,\infty)$. There can be at most one root $y\ge 1$. We may
check
\begin{eqnarray}
\Psi(x,1)&=&(x-1)\log\frac{x-1}2-x\log x,\\
\Psi(x,y)&\sim&\log y+O(1)\quad\mbox{as } y\to\infty.
\end{eqnarray}
Since $\Psi(x,1)<0$, and $\Psi(x,y)\to\infty$ as $t\to\infty$, there exists
a root $y\ge 1$ indeed.

A straightforward attempt to solve $\Psi(x,y)=0$ asymptotically gives
(\ref{gasympt}).

To prove the inequalities in (\ref{ineq4big}), we show
\begin{equation}\label{psiineq}
\Psi\big(x,\tau\,(x-1)+1\big)<0,\qquad
\Psi\left(x,\tau\,(x-1)+\frac{\tau-1}{2}\right)>0.
\end{equation}
Then the monotonicity of $y(x)$ will imply (\ref{ineq4big}).

First we show the second inequality. We substitute $y=\tau\,x-(\tau+1)/2$
into $\Psi(x,y)$:
\begin{eqnarray}
\Psi(x,y)&\!=\!&\log(\tau) +
x\,\log\left(1-\frac{\tau+1}{2\,\tau\,x}\right)-
(x-1)\,\log\left(\frac{x}{x-1}\left(1-\frac{\tau-1}{2\,\tau\,x}\right)\right)\nonumber\\
 \label{gseries1} &\!=\!&\sum_{j=1}^{\infty} \frac{1}{j\,(j+1)} \left( 1 -j
\left(\frac{\tau+1}{2\,\tau}\right)^{j+1}
-\frac{\tau\,j+j+2\,\tau}{2\,\tau}\,\left(\frac{\tau-1}{2\,\tau}\right)^j\,
\right)\frac{1}{x^j}.
\end{eqnarray}
The power series converges for $x>1$, since a tail of it can be majorated by
$\sum_{j}\frac{1}{j\,(j+1)}\,x^{-j}$. The series terms are positive for
large enough $j$.  The first terms of
(\ref{gseries1}) are 
\[
\frac{0.03017\ldots}{x^2}+\frac{0.03017\ldots}{x^3}+
\frac{0.02564\ldots}{x^4}+\ldots
\]
After applying Lemma \ref{seqlemma} twice with $p=(\tau\pm 1)/2\tau$, we
conclude that all terms in the series are positive. Hence the second
inequality in (\ref{psiineq}) follows.

If $y=\tau\,x-\tau+1$, then
\begin{eqnarray}
\Psi(x,y)&\!=\!&\log(\tau) + x\,\log\left(1-\frac{\tau-1}{\tau\,x}\right)-
(x-1)\,\log\left(\frac{x}{x-1}\left(1-\frac{\tau-2}{\tau\,x}\right)\right)\nonumber\\
 \label{gseries2} &=&\sum_{j=1}^{\infty} \frac{1}{j\,(j+1)}\,\left(1
-j\,\left(\frac{\tau-1}{\tau}\right)^{j+1}
-\frac{2\,j+\tau}{\tau}\,\left(\frac{\tau-2}{\tau}
\right)^j\,\right)\,\frac{1}{x^j}.
\end{eqnarray}
The power series converges for $x>1$, just as (\ref{gseries1}). The series
terms are positive for large enough $j$. The first terms of
(\ref{gseries2}) are 
\[
-\frac{0.10522\ldots}{x}-\frac{0.02770\ldots}{x^2}-\frac{0.00378\ldots}{x^3}+
\frac{0.00466\ldots}{x^4}+\ldots
\]
Applying lemma \ref{seqlemma} twice with $p=(\tau-1)/\tau$ and
$p=(\tau-2)/\tau$ 
we conclude that starting with the power $x^{-4}$ the coefficients are
positive. Hence the first three terms in (\ref{gseries2}) are negative, and
all remaining terms in (\ref{gseries2}) are positive. Let us consider the
function
\begin{equation}
\Psi_1(x)=x^3\,\Psi(x,\,-\tau\,x+\tau+1).
\end{equation}
The Laurent series of the derivative of this function at $x=\infty$ is:
\[
\frac{d\,\Psi_1(x)}{dx}= 
-0.21044\ldots{x}-0.02770\ldots-\frac{0.00466\ldots}{x^2}+\ldots
\]
The information about the signs of the coefficients in (\ref{gseries2})
implies that all nonzero terms in the Laurent series are negative.
Therefore $\Psi_1(x)$ is a decreasing function on the interval $(1,\infty)$.
Further, $\lim_{x\to 1^+} \Psi_1(x)=0$, since $\Psi(x,y)$ is continuous and
$\Psi(1,1)=0$. Therefore $\Psi_1(x)<0$ for $x\in(1,\infty)$. Consequently,
the first inequality in (\ref{psiineq}) follows as well.\QED

\noindent The main result about the function $\widehat{g}(t)$ defined in
(\ref{simplegt}) is the following.
\begin{theorem} \label{gttheorem}
Suppose that  $\alpha, \gamma$ are integers. 
Then
\begin{equation} \label{sigmaform}
\sup_{t>0} \widehat{g}(t)=\left\{ \begin{array}{rl}
1, & \mbox{if } \alpha=\gamma, \\
\infty, & \mbox{if } \gamma=0,\ \alpha\neq 0,\\
\max(|\frac{\alpha}{\gamma}|,\,1), & \mbox{if } \gamma>0, \\
y(\frac{\alpha}{\gamma}), & \mbox{if } \alpha<\gamma<0, \\
1+1/y(\frac{\gamma}{\gamma-\alpha}), & \mbox{if } \gamma<\alpha<0, \\
2, & \mbox{if } \gamma<0,\ \alpha=0, \\
1+y(\frac{\gamma-\alpha}{\gamma}), & \mbox{if } \gamma<0<\alpha.\\
\end{array} \right.
\end{equation}
where the function $y(x)$ is defined in Lemma $\ref{ytlemma}$.
\end{theorem}
\proof If $\alpha=\gamma$, then $\widehat{g}(t)\equiv 1$. If $\gamma\neq 0$,
then
\[ 
\widehat{g}(t)=\left|1+\frac1{\alpha\,t}\right|^{\alpha\,t}\, |\alpha\,t+1|,
\] 
and $\widehat{g}(t)\sim\exp(1)\,|\alpha|\,t$ as $t\to\infty$. 
If $\gamma>0$, we apply part \refi{iii} of Lemma \ref{gtsplemma}.

From now on we assume $\gamma<0$, $\alpha\neq\gamma$. We use Lemma
\ref{gtsplemma} with the flipped signs of $\gamma$, $t$ and $\alpha$. By
part \refi{iv}, the supremum is a local extremum, so it is achieved for some
$t=t_{\sup}$ (dependent on $\alpha$ and $\gamma$) satisfying
$\widehat{g}{\,}'(t_{\sup})=0$. Expression (\ref{ghtderiv}) gives the
following equation for $t_{\sup}$:
\begin{equation} \label{tsuprel}
\frac{|\alpha\,t_{\sup}+1|^{\alpha}\;|\gamma\,t_{\sup}|^{\gamma}}
{|\alpha\,t_{\sup}|^{\alpha}\;|\gamma\,t_{\sup}+1|^{\gamma}}=1.
\end{equation}
Hence,
\begin{equation} \label{gtsupv}
\widehat{g}(t_{\sup}) =  \frac{|\alpha\,t_{\sup}+1|}{|\gamma\,t_{\sup}+1|}.
\end{equation}
Let us define the function
\begin{equation} \label{gtsupv2}
\wtilde{y}(\alpha,\gamma)=\frac{\alpha\,t_{\sup}+1}{\gamma\,t_{\sup}+1},
\end{equation}
so that $\widehat{g}(t_{\sup})=|\wtilde{y}(\alpha,\gamma)|$. We have:
\begin{eqnarray} \label{thebluf}
\frac{|\wtilde{y}(\alpha,\gamma)|^{\alpha}}
{|\wtilde{y}(\alpha,\gamma)-1|^{\alpha-\gamma}} & = &
\frac{|\alpha\,t_{\sup}+1|^{\alpha}\;|t_{\sup}|^{\gamma-\alpha}}
{|\gamma\,t_{\sup}+1|^{\gamma}\;|\alpha-\gamma|^{\alpha-\gamma}}
\nonumber \\ & = & \frac{|\alpha|^{\alpha}}
{|\gamma|^{\gamma}\;|\alpha-\gamma|^{\alpha-\gamma}},
\end{eqnarray}
where the second equality holds because of (\ref{tsuprel}). Formula
(\ref{thebluf}) implies that $\wtilde{y}(\alpha,\gamma)$ is a real solution
of
\begin{equation} \label{eqfory}
\frac{\left|\widetilde{y}\right|^x}{|x|^x} =
\frac{\left|\wtilde{y}-1\right|^{x-1}}{|x-1|^{x-1}}, \qquad \mbox{where}
\quad x=\frac{\alpha}{\gamma}, \quad \wtilde{y}=\wtilde{y}(\alpha,\gamma)
\end{equation}
Conversely, if (\ref{eqfory}) holds, then expression (\ref{gtsupv2}) is also
true provided that $t_{\sup}$ is well defined, which is not the case only
when $\wtilde{y}(\alpha,\gamma)=\alpha/\gamma$. It follows that all
solutions of (\ref{eqfory}) except $\widetilde{y}=x$ correspond to local
extremuma of $\widehat{g}(t)$. 
We need a solution of (\ref{eqfory}) whose absolute value is greater than
$\max(1,|x|)$. 

The cases $x=1$ and $x=0$ can be proved by solving the equation
(\ref{eqfory}) directly.
\begin{figure}
\centering \epsfig{file=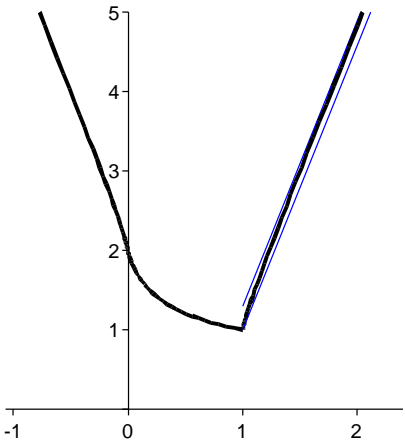, width=2.6in} \caption{The supremum of
$\widehat{g}(t)$ for $\gamma\,t<0$, as a function of
$\displaystyle\,\frac{\alpha}{\gamma}$.}
 \label{fig:cycle}
\end{figure}

If $x>1$, we have two possibilities: either $\widetilde{y}>1$ or
$\wtilde{y}<-1$ for the relevant solution of (\ref{eqfory}). But if
$\wtilde{y}>0$ and $\wtilde{y}>x$, then the left-hand side of (\ref{eqfory})
is always bigger than the right-hand side. Hence the relevant solution has
$\wtilde{y}<-1$. Then $y=|\wtilde{y}|$ satisfies (\ref{eq4:ybigx}), so the
supremum for $\alpha<\gamma<0$ is equal to $y\left(\alpha/\gamma\right)$.

If $x<0$ then the transformation $x\mapsto 1-x$, $\widetilde y\mapsto
1-\widetilde{y}$ transforms equation (\ref{eqfory}) to the same equation
with $x>1$. Since $\widetilde{y}=-y\left(\alpha/\gamma\right)$ for $x>1$, we
get the result for $\gamma<0<\alpha$.

Similarly, if $x\in(0,1)$ then the transformation $x\mapsto 1/(1-x)$,
$\widetilde y\mapsto 1/(1-\widetilde{y})$ transforms equation (\ref{eqfory})
to the same equation with $x>1$. The inverse transformation on
$\widetilde{y}$ for $x>1$ is $1-1/\widetilde{y}$, with
$\widetilde{y}=-y\left(\alpha/\gamma\right)$ again. Hence the remaining case
$\gamma<\alpha<0$ follows. \QED

Figure \ref{fig:cycle} gives the graph of $\sup_{t>0}\widehat{g}(t)$ for
$\lambda<0$ as a function of $\alpha/\gamma$, as specified by Theorem
\ref{gttheorem}. The continuous graph is piecewise defined on the intervals
$(-\infty,0)$, $(0,1)$ and $(1,\infty)$. On the interval $(1,\infty)$, the
function is identical to the function $y(x)$ of Lemma \ref{ytlemma}. The
thin lines above the interval $[1,\infty)$ are the bounding lines in (\ref{ineq4big}). 
As we see, the function approaches the asymptotic straight line very fast.
The function can be transformed between the three intervals by the fractional-linear
transformations implied in Theorem \ref{gttheorem}. The tangent slopes at
$\alpha/\gamma=1$ (from the right) and at $\alpha/\gamma=0$ (from both
sides) are actually vertical. To see this at $\alpha/\gamma=1$, compute
$dy/dx$ from $\Psi(x,y)=0$ as in (\ref{eq:psixy}). The tangent slope at
$\alpha/\gamma=1$ from the left is equal to $-1/\tau$.

As we see, the graph in Figure \ref{fig:cycle} grows rather fast with
$|\alpha/\gamma|$. If one tries to estimate the supremum of $g(t)$ by
expressing it as a product of $\widehat{g}(t)$'s, the negative $\gamma_j$'s
should be preferably paired with negative $\alpha_j$'s of similar magnitude,
so that the respective quotients $\alpha/\gamma$ would be close to 1.

\bibliographystyle{alpha}

\end{document}